\documentclass[a4paper,10pt,twoside]{article}
\usepackage{pst-all}
\usepackage{graphicx}
\usepackage{fancyhdr}
\usepackage{amsmath}
\usepackage{epsfig} 
\usepackage{psfrag} 
\usepackage{amsmath,amsfonts,amssymb,amsthm} 
\usepackage{bbm} 
\usepackage{lastpage}       

\newtheorem{Theorem}{Theorem}[section] 
\newtheorem{Corollary}[Theorem]{Corollary} 
\newtheorem{Lemma}[Theorem]{Lemma}

\newtheorem{Remark}[Theorem]{Remark}

\theoremstyle{definition} 

\theoremstyle{remark}

\newcommand{\Path}{P_{0,n-\mbox{\scriptsize path}}^p(D)}
\newcommand{\Pathst}{P_{s,t-\mbox{\scriptsize path}}^p(D)}
\newcommand{\Pathleq}{P_{0,n-\mbox{\scriptsize path}}^{\leq p}(D)}

\newcommand{\Pathleqst}{P_{s,t-\mbox{\scriptsize path}}^{\leq p}(D)}
\newcommand{\Pathgeqst}{P_{s,t-\mbox{\scriptsize path}}^{\geq p}(D)}

\newcommand{\OnePath}{P_{0,n-\mbox{\scriptsize path}}^1(D)}
\newcommand{\TwoPath}{P_{0,n-\mbox{\scriptsize path}}^2(D)}
\newcommand{\ThreePath}{P_{0,n-\mbox{\scriptsize path}}^3(D)}

\newcommand{\paths}{$(0,n)-p$-paths }
\newcommand{\fourpaths}{$(0,n)-4$-paths }

\newcommand{\UCycle}{P_C^p(K_n)}
\newcommand{\UCycleleq}{P_C^{\leq p}(K_n)}
\newcommand{\UCyclegeq}{P_C^{\geq p}(K_n)}

\newcommand{\Cycle}{P_C^p(D_n)_{| n}}
\newcommand{\Cycleo}{P_C^p(D_n)}
\newcommand{\Cycleleq}{P_C^{\leq p}(D_n)}
\newcommand{\Cyclegeq}{P_C^{\geq p}(D_n)}
\newcommand{\oez}{\"Ozl\"uk }
\newcommand{\bid}{\mbox{\textnormal{bid}} }

\newcommand{\pbowties}{$p$-bowties }
\newcommand{\pbowtie}{$p$-bowtie }
\newcommand{\UPath}{P^p_{[0,n] \mbox{ - \scriptsize{path}}}(K_{n+1}) }

\newcommand{\UPathst}{P^p_{[s,t] \mbox{ - \scriptsize{path}}}(G) }
\newcommand{\UPathleqst}{P^{\leq p}_{[s,t] \mbox{ - \scriptsize{path}}}(G) }
\newcommand{\UPathgeqst}{P^{\geq p}_{[s,t] \mbox{ - \scriptsize{path}}}(G) }
\newcommand{\UPathone}{P^1_{[0,n]  \mbox{-\scriptsize{path}}}(K_{n+1}) }
\newcommand{\UPathtwo}{P^2_{[0,n] \mbox{-\scriptsize{path}}}(K_{n+1}) }
\newcommand{\UPaththree}{P^3_{[0,n] \mbox{-\scriptsize{path}}}(K_{n+1}) }
\newcommand{\UPathn}{P^n_{[0,n] \mbox{-\scriptsize{path}}}(K_{n+1})}
\newcommand{\UPathff}{P^4_{[0,6] \mbox{-\scriptsize{path}}}(K_{6}) }
\newcommand{\upath}{$[0,n]-p$ - path }
\newcommand{\symmetric}{pseudo-symmetric }
\newcommand{\rb}[1]{\raisebox{1.5ex}[-1.5ex]{#1}}
\newcommand{\path}{$(0,n)-p$-path }
\newcommand{\pathst}{$(s,t)-p$-path }

\title{Facets of the $(s,t)-p$-path polytope}
\author{R\"udiger Stephan}
\date{}

\begin{document}

\maketitle

\begin{abstract}
\noindent
We give a partial description of the $(s,t)-p$-path polytope of a directed 
graph $D$ which is the convex hull of the incidence vectors of simple
directed $(s,t)$-paths in $D$ of length $p$. First, we point out
how the $(s,t)-p$-path polytope is located in the family of path
and cycle polyhedra. Next, we give some classes of valid inequalities
which are very similar to inequalities which are valid
for the $p$-cycle polytope, that is, the convex hull of the incidence vectors
of simple cycles of length $p$ in $D$. We give necessary and
sufficient conditions for these inequalities to be facet
defining. Furthermore, we consider a class of inequalities that has
been identified to be valid for $(s,t)$-paths of cardinality at most $p$. 
Finally, we transfer the results to related polytopes, in particular,
the undirected counterpart of the $(s,t)-p$-path polytope.
\end{abstract}

\section{Introduction}
Given a directed graph $D=(V,A)$, we say that a subset 
\[P=\{(i_1,i_2),(i_2,i_3),\dots,(i_{k-1,i_k})\}\] 
of $A$
is a directed simple $(s,t)$-path if
$k \geq 2$, all nodes $i_1,\dots,i_k$ are distinct, $s$ is the origin,
and $t$ is the terminus, that is, $s=i_1$, $t=i_k$. Below a directed
simple path will be sometimes denoted by a tuple of nodes. For
example, $(i_1,i_2,i_3,i_4)$ denotes the path
$\{(i_1,i_2),(i_2,i_3),(i_3,i_4)\}$. 
In this paper we study the facial structure of the $(s,t)-p$-path
polytope $\Pathst$ which is the convex hull of the incidence vectors
of directed $(s,t)$-paths with exactly $p$ arcs. The corresponding
\pathst problem, that is, the problem of finding a minimum cost \pathst,
is NP-hard, since for $p=n$ and negative arc cost it is equivalent to
the Hamiltonian $(s,t)$-path problem. So for general $p$ we cannot
expect to obtain a complete and tractable linear characterization of
the \pathst polytope $\Pathst$. 

A lot of path and cycle polyhedra are well studied.
Dahl and Gouveia \cite{Dahl} gave some valid inequalities for
polyhedra associated with the directed hop-constrained
shortest path problem which is the problem of finding a minimum 
$(s,t)$-path with at most $p$ arcs. Dahl and Realfsen \cite{DR}
studied the same problem on acyclic directed graphs, in particular, on
\emph{2-graphs}. The dominant of the directed 
$(s,t)$-path polytope which is the Minkowski sum of the convex hull of the
incidence vectors of simple $(s,t)$-paths and the Euclidean space
$\mathbb{R}^A$ is determined by nonnegativity constraints $x_{ij} \geq
0$ and cut inequalities $x(C) \geq 1$ for all 
$(s,t)$-cuts $C$ (see Schrijver \cite{Schrijver2003}, chapter 13). The
cycle polytope $P_C(D_n)$ which is the convex hull of the incidence
vectors of all simple directed cycles of the complete directed graph
$D_n$ has been investigated by Balas and Oosten \cite{BO}, while the
undirected counterpart, the circuit polytope, has been studied by
Coullard \& Pulleyblank \cite{CP} and Bauer \cite{Bauer}. Hartmann and
\oez \cite{HO} gave a partial description of the $p$-cycle polytope
$P_C^p(D_n)$ which is the convex hull of the incidence vectors of all
simple $p$-cycles of $D_n$. 
Maurras and Nguyen \cite{MN1,MN2} studied the facial
structure of the undirected analog. Finally, Bauer et al. \cite{BLS}
also studied the cardinality constrained circuit polytope, which is
the convex hull of the incidence vectors of all undirected simple
cycles with at most $p$ edges on the complete graph $K_n$.

The present paper is motivated by the observation that the \pathst
polytope is closely related to the polyhedra mentioned in the last paragraph
and it has an exposed position among them. Indeed, valid inequalities
for the \pathst polytope can easily be transformed into valid
inequalities for some related polytopes, for example, by lifting. A first
overview is given in Figure 1. An arrow there
means that facet defining inequalities (or some classes of facet
defining inequalities) of the polytope at the tail of the arrow can be
transformed into facet defining inequalities for the polytope at the
head of the arrow, where $G$ and $D$ are appropriate graphs and
digraphs, respectively.

The remainder of the paper is organized as follows: In Section 2 we
propose an integer programming formulation of the \pathst polytope and
describe how valid inequalities can be lifted to valid inequalities of
the $p$-cycle polytope. Section 3 contains the study of the facial
structure of the \pathst polytope $\Pathst$ on an appropriate digraph
$D$. Finally, in Section 4 we transfer the results of Section 3 to the
polytopes mentioned that are related to the \pathst polytope $\Pathst$. 

\begin{center}
\psset{xunit=0.7cm,yunit=1cm,runit=1cm} 
\begin{pspicture}(0,0)(16,7) 
\rput(0,0){\rnode{1}{ \psframebox{$\UPathleqst$}}} 
\rput(3,1){\rnode{2}{ \psframebox{$\UPathst$}}} 
\rput(6,2){\rnode{3}{ \psframebox{$\UPathgeqst$}}} 

\rput(0,4){\rnode{4}{ \psframebox{$\UCycleleq$}}} 
\rput(3,5){\rnode{5}{ \psframebox{$\UCycle$}}} 
\rput(6,6){\rnode{6}{ \psframebox{$\UCyclegeq$}}} 

\rput(10,0){\rnode{7}{ \psframebox{$\Pathleqst$}}} 
\rput(13,1){\rnode{8}{ \psframebox*[fillcolor=yellow]{$\Pathst$}}} 
\rput(13,1){\rnode{8}{ \psframebox{$\Pathst$}}} 
 
\rput(16,2){\rnode{9}{ \psframebox{$\Pathgeqst$}}} 

\rput(10,4){\rnode{10}{ \psframebox{$\Cycleleq$}}} 
\rput(13,5){\rnode{11}{ \psframebox{$\Cycleo$}}} 
\rput(16,6){\rnode{12}{ \psframebox{$\Cyclegeq$}}} 

\ncline{->}{2}{1}
\ncline{->}{2}{3}
\ncline[linestyle=dashed]{->}{2}{5} \naput{?}
\ncline{->}{5}{4}
\ncline{->}{5}{6}
\ncline{->}{8}{7}
\ncline{->}{8}{9}
\ncline{->}{8}{11}
\ncline{->}{11}{10}
\ncline{->}{11}{12}
\ncline{->}{11}{5} 
\ncline{->}{8}{2}  

\rput(8,-1.5){Figure 1. The \pathst polytope $\Pathst$ and related polytopes.}
\end{pspicture}
\end{center}

\newpage

\section{Basic results}

We start the polyhedral analysis of the \pathst polytope with an
integer programming formulation. In the sequel, $D=(V,A)$ is a digraph
on node set $V=\{0,\dots,n\}$  whose arc set $A$ contains neither
loops nor parallel arcs. The nodes $s$ and $t$ will be
identified with the nodes $0$ and $n$, respectively. Consequently, the
\path polytope will be denoted by $\Path$. The integer points 
of $\Path$ are characterized by the system 
\begin{align}
x(\delta^-(0)) & = 0, & \label{zero1}  \\ 
x(\delta^+(n)) & =0, &  \label{zero2} \\
x(\delta^+(i))- x(\delta^-(i)) & = 
 \left \{ \begin{array}{r@{}l}
1 & \mbox{ if } i=0,\\
0 &\mbox{ if }  i \in V \setminus \{0,n\},\\
-1 & \mbox{ if }  i=n,\\
\end{array} \right. \label{flow}\\
x(A) &=p, \label{card}\\
x(\delta^+(i)) & \leq 1 & \forall \,\,\,i \in V\setminus \{0,n\}, \label{degree}\\
x((S:V \setminus S)) & \geq x(\delta^+(j)) & \forall  S \subset
V, 3 \leq |S| \leq n-2, \label{osmincut}\\
\nonumber & &  0,n \in S, j \in V \setminus S,\\
x_{ij} & \in \{0,1\}  & \forall \,\,\, (i,j) \in A. \label{integer}
\end{align}
Here, we denote by $\delta^+(k)$ and $\delta^-(k)$  the set of arcs directed
out of and into node $k$, respectively. For an arc set $F \subseteq A$
we set $x(F):=\sum_{(i,j) \in F} x_{ij}$, and for any node sets $S,T$ of
$V$, $(S:T)$ is short for $\{(i,j) \in A| i \in S, j \in
T\}$. Furthermore, in the following we denote by $A(S)$ the subset of arcs
whose both endnodes are in $S$, for some $S \subseteq V$. 

The incidence vectors of node-disjoint unions of a $(0,n)$-path and
cycles on node set $V \setminus \{(0,n)\}$ are described by the
equations (\ref{zero1})-(\ref{zero2}), the
\emph{flow constraints} (\ref{flow}), \emph{degree constraints}
(\ref{degree}), and the
\emph{integrality constraints} (\ref{integer}). The \emph{one-sided
  min-cut inequalities} (\ref{osmincut}) are satisfied by all
$(0,n)$-paths but violated by unions of a $(0,n)$-paths and cycles
on $V \setminus \{0,n\}$. Finally, the \emph{cardinality constraint}
(\ref{card}) ensures that all $(0,n)$-paths are of length $p$.

Complete linear descriptions of $\Path$ for $p=1,2,3$ are given in
Table 1, where $D$ is the complete digraph on node set
$\{0,\dots,n\}$. The results for $p=2$ and $p=3$ follows from the fact
that a $(0,n)-2$-path visits exactly one internal node and a
$(0,n)-3$-path contains exactly one internal arc. Since the number of
internal nodes is $n-1$, the dimension of $\TwoPath$ is $n-2$, and
since the number of internal arcs is $(n-1)(n-2)$, the dimension of
$\ThreePath$ is equal to $(n-1)(n-2)-1=n-3n+1$. 
The $(0,n)-1$-path polytope $\OnePath$ has clearly dimension $0$ and
is determined by the equations $x_{0n}=1$ and $x_{ij}=0$ for all
$(i,j) \in A \setminus \{(0,n)\}$. We suppose in the sequel that $A$
contains all arcs $(i,j)$, where $i \neq j \in V$, except the arcs
$(i,0)$, $(n,i)$ for $i=1,\dots,n-1$, $(0,n)$, and $(n,0)$. 

Contracting the nodes $0$ and $n$ to the single node $n$ we obtain the
complete digraph $D_n$ on $n$ nodes, and we see that the set of simple
$(0,n)-p$-paths defined on $D$ can be identified with the set of simple
$p$-cycles defined on $D_n$ that contain node $n$. Hence, 
the $(0,n)-p$-path polytope $\Path$ and the node constraint cycle
polytope $\Cycle:=\{\textbf{x} \in P_C^p(D_n)| x(\delta^+(n))=1\}$ are
isomorphic. 
In particular, when $p=n$, $\Path$ is isomorphic to the asymmetric
traveling salesman polytope which has
dimension $n^2-3n+1$ (see \cite{Groetschel}). Furthermore, Hartmann and
\"Ozl\"uk \cite{HO} showed that $\Cycle$ is a facet of the
$p$-cycle polytope if $4 \leq  p <n$. For $4 \leq p <n$, the $p$-cycle
polytope has dimension $n^2-2n$ and therefore the dimension of $\Path$
is equal to $n^2-2n-1$. Moreover, this relation leads to the following
theorem obtained by standard lifting (see Nemhauser and Wolsey \cite{NW}).

\begin{table}
\noindent Table 1.  Polyhedral Analysis of $\Path$, where $D$ is the complete
digrah on node set $\{0,\dots,n\}$. For $p=n$, $\Path$ is equivalent
to the asymmetric traveling salesman polytope defined on $n$ nodes.

\vspace{0.4 cm}
\noindent \begin{tabular}{|c|c|rll|} \hline
& & & & \\ 
\rb{$p$} & \rb{Dimension} & \multicolumn{3}{c|}{\rb{Complete linear description}} \\ \hline 
      &  & $x_{0n}$ & $=1$ & \\
\rb{$1$} & \rb{0} & $x_{ij}$ & $=0$ & $\forall \; (i,j) \in A \setminus \{(0,1)\}$ \\ \hline 
   &  & $x(\delta^-(0))$ & $=0$ &  \\
      &       & $x(\delta^+(n))$ & $=0$ & \\
   &    & $x_{ij}$   & $=0$ & $\forall \: (i,j) \in A(V \setminus \{0,n\})$ \\
\rb{$2$} & \rb{$n-2$} & $x(\delta^+(0))$ & $ =1$ & \\ 
   &    & $x_{0j}-x_{jn}$ & $=0$ & $\forall \: j \in V \setminus \{0,n\}$\\
& & $x_{0j}$   & $\geq 0$ & $\forall \: j \in V \setminus \{0,n\}$ \\ \hline
    &  & $x(\delta^-(0))$ & $=0$ & \\
      &       & $x(\delta^+(n))$ & $=0$ & \\
& & $x(A(V \setminus \{0,n\}))$ &$=1$ & \\
\rb{$3$}& \rb{$n^2-3n+1$} & $x(\delta^+(i))$ & $= x_{0i}+x_{in}$ & $ \forall \: i \in V \setminus \{0,n\}$\\
& & $x(\delta^-(i))$ & $=x_{0i}+x_{in}$ & $\forall \: i \in V \setminus \{0,n\}$\\
& & $x_{ij}$ & $\geq 0$  & $\forall \: (i,j) \in A(V \setminus \{0,n\})$ \\  \hline \hline
& & & & \\
 &  & \multicolumn{3}{c|}{\rb{Partial linear description}} \\ \hline
$4$ & & & & \\
$\vdots$ & $n^2-2n-1$ & \multicolumn{3}{c|}{ \rb{equations (\ref{zero1})-(\ref{card})}} \\
$n-1$ & & \multicolumn{3}{c|}{\rb{see Section 3}} \\ \hline
\end{tabular}
\end{table}

\begin{Theorem} \label{T0}
Let $\textbf{ax} \leq a_0$ be a facet defining inequality for the $(0,n)-p$-path polytope $\Path$, where $4 \leq p< n$, and let $\gamma$ be the maximum
of $a(C)$ over all $p$-cycles $C$ in $D$. Setting $a_{ni}:=a_{0i}$ for $i=1,\dots,n-1$, the inequality
\begin{equation}
\sum_{i=1}^n \sum_{j=1 \atop j \neq i}^n a_{ij} x_{ij} + (\gamma -a_0)x(\delta^+(n)) \leq \gamma
\end{equation}
defines a facet of the $p$-cycle polytope $P_C^p(D_n)$, where $D_n$ is the complete digraph on node set $\{1,\dots,n\}$. \hfill $\Box$
\end{Theorem}

This easy but fundamental relation between the \path polytope
$\Path$ and the $p$-cycle polytope $P_C^p(D_n)$ also holds between
other length restricted path and cycle polytopes (see
\cite{Stephan}). This fact implies that it would be profitably to
study first the facial structure of a length restricted directed path
polytope and afterwards that of the corresponding cycle polytope. 
In our special case, the $p$-cycle polytope is already well studied;
so we will proceed in the opposite direction, that is, starting from
the results for the $p$-cycle polytope $P_C^p(D_n)$ given by Hartmann and \oez
\cite{HO} we will prove in many cases analogous results for the \path
polytope $\Path$ and it is not surprising that this can be often done along the
lines of the proofs of the authors mentioned above. Lemma \ref{L2}
adapts Lemmas 2 and 6 of Hartmann and \oez \cite{HO} for our purposes.
The other statements of this section can be proved in the same manner
as the original statements in \cite{HO}; so we omit their proofs.

\begin{Lemma}[cf. Lemmas 2 and 6 of Hartmann and \oez \cite{HO}] \label{L2} 
Let $3 \leq p < n$, $\textbf{c}$ be a row vector, $s,t \in V$, $s \neq
t$, and $R \subseteq V \setminus \{s,t,0,n\}$.
There are $\lambda$, $\pi_s$, $\pi_t$, and $\{\pi_j|j \in R\}$ with
\[\begin{array}{rcll}
c_{si} & = & \lambda + \pi_s - \pi_i & \forall \: i \in R,\\
c_{it} & = & \lambda + \pi_i - \pi_t & \forall \: i \in R,\\
c_{ij} & = & \lambda + \pi_i - \pi_j & \forall \: (i,j) \in A(R),
\end{array}
\]
if one of the following conditions holds:
\begin{itemize}
\item[(i)] $|R| \geq 5$ and $c_{ik}+c_{kj}=c_{il}+c_{lj}$ for all
  distinct nodes $i \in R \cup \{s\}$, $j \in R \cup \{t\}$, $k,l \in R$. \label{L2i}
\item[(ii)] $|R| \geq p \geq 4$ and $c(P)= \gamma$ for all $(s,t)-p$-paths
  $P$, whose internal nodes are all in $R$.  \label{L2ii}
\item[(iii)] $|R|=p-1$, $c(P)= \gamma$ for all $(s,t)-p$-paths $P$, whose internal nodes are all the nodes of $R$, and
$c(P)=\delta$ for all $(s,t)-r$-paths $P$, all $r-1$ of whose internal are in $R$, for some $2 \leq r <p$.  \label{L2iii}
\item[(iv)] $p=3$, $|R| \geq 3$, $c(P)=\gamma$ for all $(s,t)-3$-paths $P$, whose internal nodes are all in $R$, and
$c(P)=\delta$ for each $(s,t)-2$-path $P$ whose inner node is in $R$.  \label{L2iv}
\end{itemize}
\end{Lemma}

\begin{proof}
(i) In particular, $c_{ik}+c_{kj}=c_{il}+c_{lj}$ for all distinct nodes $i,j,k,l \in R$. Using Lemma 2 of Hartmann and \oez \cite{HO}, it follows that there 
are $\lambda$ and $\{\pi_j|j \in R\}$ with
\[\begin{array}{rcll}
c_{ij} & = & \lambda + \pi_i -\pi_j & \forall \: (i,j) \in A(R).
\end{array}
\]
Next, setting $\pi_s:=c_{sk}+\pi_k-\lambda$ and $\pi_t:=\lambda+\pi_k-c_{kt}$ for some $k \in R$,
we derive 
\begin{eqnarray}
\nonumber c_{si} & = & c_{sk}+c_{kl}-c_{il} = \lambda + \pi_s - \pi_i, \\
\nonumber c_{it} & = & c_{kt}+c_{lk}- c_{li} = \lambda +\pi_i -\pi_t
\end{eqnarray}
for all $i \in R$.

\vspace{0.25cm}
(ii) First, let $|R| \geq 5$. Since $|R| \geq p$, for all distinct nodes $i,j,k,l \in R$ there is a $(s,t)-p$-path that contains the arcs $(i,k)$ and $(k,j)$ but does not visit node $l$. Replacing node $k$ by node $l$ in $P$ yields another $(s,t)-p$-path and thus
$c_{ik}+c_{kj}=c_{il}+c_{lj}$ for all distinct nodes $i,j,k,l \in R$. Lemma 2 of Hartmann and \oez implies that there are
$\lambda$ and $\{\pi_j|j \in R\}$ such that $c_{ij}=\lambda +\pi_i-\pi_j$ for all $(i,j) \in A(R)$.
Set $\pi_s:=c_{sk}+\pi_k-\lambda$ and $\pi_t:=\lambda+\pi_l-c_{lt}$ for some $k \neq l \in R$.
Any $(s,t)-p$-path whose internal nodes are in $R$ and that uses the arcs $(s,k), (l,t)$ yields $\gamma=p \lambda +\pi_s-\pi_t$.
Further, considering for $i \in R$ a $(s,t)-p$-path $P$ whose internal nodes are in $R$ and that uses the arcs $(s,i), (l,t)$
yields $c_{si}=\lambda+\pi_s-\pi_i$ for all $i \in R$. Analogous it follows that $c_{jt}=\lambda+\pi_j-\pi_t$ for all $j \in R$.

Next, let $|R|=p=4$. Without loss of generality, we may assume that $R= \{1,2,3,4\}$. Setting $Q:= \{1,2,3\}$ and identifying the nodes $s$ and $t$, Theorem 23 of Gr\"otschel and Padberg implies  that there are $\alpha_s$, $\beta_t$, $\{\alpha_j| j \in Q\}$, and $\{\beta_j|j \in Q\}$ such that
\[\begin{array}{rclcl}
c_{si} & = & \alpha_s + \beta_i & & \forall \: i \in Q,\\
c_{ij} & = & \alpha_i +\beta_j & & \forall \: (i,j) \in A(Q),\\
c_{it} & = & \alpha_i +\beta_t && \forall \: i \in Q.
\end{array}
\]
Considering for any two nodes $i \neq j \in Q$ the $(s,t)-4$-paths $(s,4,k,i,t)$ and $(s,4,k,j,t)$, where $k$ is the remaining node in $Q$, we see that $c_{ki}+c_{it}=c_{kj}+c_{jt}$ which implies that $\alpha_i+\beta_i=\alpha_j+\beta_j$ for all $i,j \in Q$. Denoting by $\lambda$ this common value and setting $\pi_s:=\alpha_s$, $\pi_j:=\alpha_j$ for $j=1,2,3$, and $\pi_t:= \lambda-\beta_t$, yields $c_{si}=\lambda+\pi_s-\pi_i$, $c_{it}=\lambda+\pi_i-\pi_t$ for $i=1,2,3$, and $c_{ij}=\lambda+\pi_i-\pi_j$ for all $(i,j) \in A(Q)$. Now setting
$\pi_4:=\lambda+\pi_s-c_{s4}$, we see that $c_{4t}=\lambda+ \pi_4-\pi_t$, $c_{i4}=\lambda+\pi_i-\pi_4$, and $c_{4i}=\lambda+\pi_4-\pi_i$ for $i=1,2,3$.

\vspace{0.25cm}
(iii) This is Lemma 6 of Hartmann and \oez \cite{HO}. 

\vspace{0.25cm}
(iv) Without loss of generality, let $1,2 \in R$. Condition (iii)
implies that there are $\lambda$, $\pi_s$, $\pi_1$, $\pi_2$, and
$\pi_t$ with the required property restricted on
$Q:=\{1,2\}$. Further, it follows that $\gamma=3 \lambda+\pi_s-\pi_t$
and $\delta=2 \lambda +\pi_s-\pi_t$. Setting
$\pi_i:=\lambda+\pi_s-c_{si}$ for all $i \in R \setminus Q$, we see
immediately that $c_{it}=\lambda+ \pi_i \pi_t$ for all $i \in R
\setminus Q$. Thus we also obtain $c_{ij}=\lambda+\pi-\pi_j$ for all
$(i,j) \in A(R)$.
\end{proof}

Equivalence of inequalities is an important matter when studying
polyhedra. Two valid inequalities for the \path polytope $\Path$ are
equivalent if one can be obtained from the other by multiplication
with a positive scalar and adding appropriate multiples of the flow
conservation constraints (\ref{flow}) and the cardinality constraint
(\ref{card}). Clearly, two valid inequalities define the same facet of $\Path$
if and only if they are equivalent. For the next theorem that 
establishes a relationship between a linear basis of equality system
(\ref{flow}), (\ref{card}) and the arcs defining it we introduce the
following two definitions: a \emph{balanced cycle} is a (not 
necessarily directed) simple cycle that contains the same number of
forward and backward arcs and an \emph{unbalanced 1-tree} is a
subgraph of $D$ consisting of a spanning tree $T$ plus an arc $(k,l)$
whose fundamental cycle $C(k,l)$ is not balanced.

\begin{Theorem}[cf. Theorem 3 of Hartmann and \oez \cite{HO}]  \label{T3}
Let $n \geq 2$ and let $H$ be a subgraph of $D$. The variables
corresponding to the arcs of $H$ form a basis for the linear equality
system \emph{(\ref{flow})}, \emph{(\ref{card})} if and only if $H$ is an unbalanced
1-tree. \hfill $\Box$ 
\end{Theorem}

\begin{Corollary}[cf. Corollary 4 of Hartmann and \oez \cite{HO}]  \label{C4}
Let $\textbf{cx} \leq c_0$ be a valid inequality for $\Path$, and let
values $b_{ij}$ be specified for the arcs $(i,j)$ in an unbalanced
1-tree $H$. Then there is an equivalent inequality $\textbf{c'x} \leq
c'_0$ for which $c'_{ij} = b_{ij}$ for all arcs $(i,j)\in H$. 
\hfill $\Box$ 
\end{Corollary}

\begin{Corollary}[cf. Corollary 5 of Hartmann and \oez \cite{HO}] \label{C5}
Let $3 \leq p <n$, $\textbf{c}$ be a row vector, $s \in V \setminus
\{n\}$, $t \in V \setminus \{0\}$, $s \neq t$, $R \subseteq V \setminus
\{s,t,0,n\}$ with $|R| \geq 2$, let either of the conditions of Lemma
\emph{\ref{L2}} be satisfied, and suppose that $c_{ij}= \beta$ holds for all
$(i,j)$ in an unbalanced 1-tree $H$ on $R$. Then $c_{ij}= \beta$ for
all $i,j \in R$. Moreover, there are $\sigma$ and $\tau$ with
$c_{si}=\sigma$ and $c_{it}=\tau$ for all $i \in R$. 
\end{Corollary}

\begin{proof}
In either case, Lemma \ref{L2} implies that there are
 $\lambda$, $\pi_s$, $\pi_t$, and $\{\pi_j|j \in R\}$ with
\[\begin{array}{rcll}
c_{si} & = & \lambda +\pi_s - \pi_i & \forall \: i \in R,\\
c_{it} & = & \lambda+\pi_i-\pi_t & \forall \: i \in R,\\
c_{ij} & = & \lambda + \pi_i -\pi_j & \forall \: (i,j) \in A(R),
\end{array}
\]
Without loss of generality, let $\pi_k = 0$ for some $k \in
R$. Theorem \ref{T3} then implies that $\lambda= \beta$ and $\pi_j=0$
for all $j \in R$. Thus, $c_{si}=\beta+ \pi_s$ and $c_{it}=\beta -
\pi_t$ for all $i \in R$. 
\end{proof}

The next theorem can be used to lift facet defining inequalities for
the \path polytope $\Path$ into facet defining inequalities for
$P_{0,n-\mbox{\scriptsize path}}^p(D')$, where $D'=D_{n+k+1} -
(\delta^-(0) \cup \delta^+(n))$. Before stating it we need some definitions.
A subset $B \subseteq A$ of cardinality $p$ is called a
\emph{p-bowtie} if it is the union of a $(0,n)$-path $P$ and a simple
cycle $C$ connected at exactly one node. The $p$-bowtie $B$ is said to be
\emph{tied} at node $k$ if $V(P) \cap V(C)=\{k\}$. 
A facet $F$ of $\Path$ is called \emph{regular} if it is defined by an
inequality $\textbf{cx} \leq c_0$ that is not equivalent to a
nonnegativity constraint $x_{ij} \geq 0$ or a \emph{broom inequality}
\begin{equation}
x((\delta^+(i)) \geq x_{ji}+x_{ik}
\end{equation}
for some internal node $i$, where $j=k$ is an internal node or $j=0$
and $k=n$. Note that $F$ is already regular if for each internal node
$k$, there is a  \path $P$ with $c(P)< c_0$ that does not visit node
$k$ (see \cite{HO}). 

\begin{Theorem}[cf. Theorem 8 of Hartmann and \oez \cite{HO}]  \label{T8}
Suppose that $\textbf{cx} \leq c_0$ induces a regular facet of $\Path$,
where $3<p<n$. Let $k$ be an internal node such that $c(B) \leq c_0$
for all $p$-bowties $B$ tied at node $k$ and let $\delta_k$ be the
maximum of $c(\Gamma)$ over all $0,n$-paths $\Gamma$ of length $p-1$
that visit node $k$. Then 
\begin{equation} \label{lift1}
\mathbf{cx} + \sum_{i=0 \atop i \neq k}^{n-1} c_{ik} x_{i,n+1}+
\sum_{j=1 \atop j \neq k}^{n}
c_{kj}x_{n+1,j} + (c_0 - \delta_k)[x_{k,n+1}+x_{n+1,k}]
\leq c_0 
\end{equation}
defines a regular facet of $P_{0,n-\mbox{\scriptsize path}}^p(D')$,
where $D'$ is the digraph obtained by subtracting from the complete
digraph on node set $\{0,\dots,n+1\}$ the arc sets $(\delta^-(0)$ and
$\delta^+(n))$.  
\hfill $\Box$
\end{Theorem}

Since inequality (\ref{lift1}) is obtained by copying the coefficient
structure of node $k$, one refers to this process as ``lifting by
cloning node $k$''. In order to show that a class $\mathcal{K}$ of
regular inequalities define facets of the \path polytope it suffices
to show it for a subclass $\mathcal{K}'\subset \mathcal{K}$ from which
the remaining inequalities in $\mathcal{K} \setminus \mathcal{K}'$ can
be obtained by cloning internal nodes. The members of a minimal
subclass $\mathcal{K}'$ (minimal with respect to set inclusion) are
said to be \emph{primitive}. 

Before stating the last theorem of this section we also need some
definitions. Let $F$ be a subset of $A$, the \emph{auxiliary graph}
$G_F$ is an undirected bipartite graph on $2n$ nodes
$v_0,\dots,v_{n-1}$, $w_{1},\dots,w_n$, with the property that $(i,j)
\in F$ if and only if $G_F$ contains the arc $(v_i,w_j)$. Given a
valid inequality $\mathbf{cx} \leq c_0$, a \path $P$ is said to be
\emph{tight} if $c(P)=c_0$. Moreover, we define the following
equivalence relation on the arc set $A$: two arcs $(i,j)$ and $(k,l)$
are related with respect to $\textbf{cx} \leq c_0$, if there is an arc
$(f,g) \in A$ with $a_{ij}=a_{fg}=a_{kl}$ and two tight
$(0,n)-p$-paths $P_{ij}$, $P_{kl}$ such that $(i,j),(f,g) \in P_{ij}$
and $(k,l), (f,g) \in P_{kl}$. 

\begin{Theorem}[cf. Theorem 9 of Hartmann and \oez \cite{HO}] \label{T9}
Let $\textbf{a} \geq \textbf{0}$ and $\textbf{ax} \leq a_0$ be a facet
defining inequality for $\Path$, where $3 < p<n$. Suppose that the
auxiliary graph $G_Z$ for the arc set $Z:=\{(i,j) \in A| a_{ij}=0\}$ is
connected, every tight \path with respect to $\textbf{ax} \leq a_0$
contains at least one arc $(i,j) \in Z$, and every arc $(i,j)$ belongs
to the same equivalence class with respect to $\textbf{ax} \leq
a_0$. Let $R$ be a set of nodes, set $q:=p+|R|$, and let $t$ be the
smallest number such that 
\begin{equation} \label{mlift}
\textbf{ax} + t \sum_{j \in R} x(\delta^+(j)) \leq a_0 +|R|t
\end{equation}
is valid for all $(0,n)-q$-paths on $V \cup R$, and if $|R| \geq 2$
suppose that at least one tight $(0,n)-q$-path with respect to
(\ref{mlift}) visits $r$ nodes in $R$ with $0<r<|R|$. Then
(\ref{mlift}) is facet defining for the $(0,n)-q$-path polytope on $V
\cup R$. 
\hfill $\Box$
\end{Theorem}

\section{Facets and valid inequalities}

In the sequel we will show that the inequalities given in the
IP-formulation, the \emph{nonnegativity constraints} $x_{ij} \geq 0$,
as well as some more inequalities are in general facet defining for
$\Path$. Throughout, we assume that $4 \leq p \leq n-1$. The
inequalities considered in 
Theorems \ref{T10} - \ref{T15} were shown to be valid for the
$p$-cycle polytope in Hartmann and \oez \cite{HO}. So they are also
valid for $\Path$, since the $(0,n)-p$-path polytope on $D$ can be
interpreted as the restriction of the $p$-cycle polytope on $D_n$ to
the hyperplane defined by $x(\delta^+(n))=1$. 

\subsection{Trivial inequalities}
\begin{Theorem}[cf. Theorem 10 of Hartmann and \oez \cite{HO}] \label{T10}
The nonnegativity constraint
\begin{equation} \label{nn}
x_{ij} \geq 0
\end{equation}
is valid for $\Path$ and induces a facet of $\Path$ whenever $4 \leq p
\leq n-1$. 
\end{Theorem}

\begin{proof}
When $n \leq 6$ and $p=4$ or $p=5$, (\ref{nn}) can be proved to induce a
facet by application of a convex hull code (e.g. Polymake
\cite{polymake}), so we assume that $n \geq 7$. Suppose that
$\textbf{cx}=c_0$ is satisfied by every $\textbf{x} \in \Path$ with
$x_{ij}=0$. At least one of the two nodes $i$ and $j$ is an internal
node, because $(0,n) \notin A$. Without loss of generality, we may
assume that $j \in \{1,\dots,n-1\}$ and set $R:=V \setminus
\{0,n,j\}$. By Corollary \ref{C4}, we may assume that
$c_{jw}=c_{0w}=c_{wn}=0$ for some $w \in R$ and $c_{kl}=0$ for all
arcs $(k,l)$ in some unbalanced 1-tree on $R$. 

Let $q \in R \cup \{0\}$, $r,s \in R$, $t \in R \cup \{n\}$ be
distinct nodes and let $P$ be a \path   that contains the arcs $(q,r)$
and $(r,t)$ but does not visit node $s$ or use the arc
$(i,j)$. Substituting node $r$ by node $s$ in $P$ we obtain another
\path that does not use $(i,j)$. Hence condition (\ref{L2i}) of Lemmma \ref{L2}
holds and Corollary \ref{C5} implies that $c_{kl}=0$ for all $(k,l)
\in A(V \setminus \{j\})$ which also implies that $c_0=0$. 

Each \path that uses the arc $(j,w)$ but does not use the arc $(i,j)$
also satisfies (\ref{nn}) with equality, so $c_{kj}=0$ for all $k \in
V \setminus \{i,n,w\}$. Similar considerations yield $c_{jk}=0$ for
all $k \in V \setminus \{0\}$ and $c_{wj}=0$ if $w \neq i$. Thus,
$c_{kl}=0$ for all arcs $(k,l) \neq (i,j)$ and therefore
$\textbf{cx}=c_0$ is simply $c_{ij}x_{ij}=0$. 
\end{proof}

\begin{Theorem}[cf. Theorem 11 of Hartmann and \oez \cite{HO}] \label{T11}
Let $j$ be an internal node. The degree constraint
\begin{equation} \label{degree1}
x(\delta^+(j)) \leq 1
\end{equation}
is valid for $\Path$ and induces a facet of $\Path$ whenever $4 \leq p
\leq n-1$. 
\end{Theorem}

\begin{proof}
Without loss of generality, we will show that $x(\delta^+(1)) \leq 1$
defines a facet of $\Path$. First we will show that Theorem \ref{T11}
holds when $p=4$. If $n=5$, $x(\delta^+(1)) \leq 1$ can be proved to
define a facet using a convex hull code.  Theorem \ref{T8} applied to
node 2 yields then the result when $n \geq 6$.

Secondly, we will investigate the case $p \geq 5$. Suppose that
$\textbf{cx}=c_0$ is satisfied by every $\textbf{x} \in \Path$ with
$x(\delta^+(1)) = 1$. By Corollary \ref{C4}, we may assume that
$c_{21}=c_{02}=c_{2n}=0$ and $c_{ij}=0$ in some unbalanced 1-tree on
$R:=\{2,3,\dots,n-1\}$. Since $|R| \geq p-1 \geq 4$ and $c(P)=c_0 -
c_{01}$ for all $(1,n)$-paths $P$ of length $p-1$ whose internal nodes
are all in $R$, condition (\ref{L2ii}) of Lemma \ref{L2} holds. Thus,
$c_{ij}=0$ for all $(i,j) \in A(R \cup \{n\})$ and $c_{1j}=1$ for all
$j \in R$ using Corollary \ref{C5}. Now it is easy to see that
$\textbf{cx}=c_0$ is simply $x(\delta^+(1))=1$. 
\end{proof}

\subsection{Cut inequalities}

\begin{Theorem}[cf. Theorem 12 of Hartmann and \oez \cite{HO}] \label{T12}
Let $S \subset V$ and $0,n \in S$. The min-cut inequality
\begin{equation} \label{mincut}
x((S:V \setminus S)) \geq 1
\end{equation}
is valid for $\Path$ if and only if $|S| \leq p$ and facet defining
for $\Path$ if and only if $3 \leq |S| \leq p$ and $|V \setminus S| \geq 2$.
\end{Theorem}

\begin{proof}
The min-cut inequality (\ref{mincut}) is valid for $\Path$ if and only
if $|S| \leq p$, since a \path can be obtained in $S$ if and only if
$|S| \geq p+1$. When $|S|=2$, (\ref{mincut}) is an implicit
equation. When  $|V \setminus S|=1$, $n \leq p$. So we suppose that $3
\leq |S| \leq p$ and $|V \setminus S| \geq 2$.

First let $|S|=3$. When $|V \setminus S| \leq 4$, (\ref{mincut}) can
be shown to be facet defining by means of a convex hull code, so let
$|V \setminus S| \geq 5$. Let w.l.o.g. $S= \{0,1,n\}$ and suppose that
$\textbf{cx}=c_0$ is satisfied by every $\textbf{x} \in \Path$ that
satisfies (\ref{mincut}) with equality. Using Corollary \ref{C4}, we
may assume that $c_{01}=0$, $c_{0w}=c_0$ and $c_{wn}=0$ for some $w
\in V \setminus S$, as well as $c_{ij}=0$ for all arcs $(i,j)$ in some
unbalanced 1-tree $H$ on $V \setminus S$.

Let $i \in (V \setminus S) \cup \{0\}$, $j \in  (V \setminus S) \cup \{n\}$, $k,l \in V \setminus S$ be distinct nodes 
and let $P$ be a tight \path that contains the arcs $(i,k),(k,j)$ but does not visit node $l$. 
Such a path $P$ exists even when $p=4$. Replacing node $k$ by node $l$ yields another tight \path, 
and hence condition (\ref{L2i}) of Lemma \ref{L2} holds. 
Corollary \ref{C5} implies that $c_{ij}=0$ for all $(i,j) \in A(V \setminus S)$, $c_{0i}=c_0$, and $c_{in}=0$ for all $i \in V \setminus S$.  
Now it is easy to see that $c_{1i}=c_0$ and $c_{i1} + c_{1n}=0$ for all $i \in V \setminus S$. Subtracting  $c_{1n}$ times the equation $x(\delta^-(n))=1$ and adding
$c_{1n}$ times the equation $x((V \setminus S:S)) - x((S: V \setminus S))=0$, we see that
$\textbf{cx}=c_0$ is equivalent to $(c_0-c_{1n}) x(S:V \setminus S)=c_0-c_{1n}$.

Secondly, let $|S| \geq 4$. 
Let w.l.o.g. $S= \{0,1,2,\dots,q,n\}$ for some $q <p$ and suppose that 
$\textbf{cx}=c_0$ is satisfied by every $\textbf{x} \in \Path$ that satisfies (\ref{mincut}) with equality.
Using Corollary \ref{C4}, we may assume that 
$c_{01}=c_{1n}=0$, $c_{1i}=c_0$ for all $i \in (V \setminus S)$, and
$c_{ij}=0$ for all arcs $(i,j)$ in some unbalanced 1-tree on $R:= S \setminus \{0,n\}$.

Let $P$ be the path $(q+1,\dots,p-1,n)$ and $Q$ be the path $(q+1,\dots,p,n)$
Then $c(\Gamma)=c_0 - c(P)$ for all $(0,q+1)$-paths $\Gamma$, whose internal nodes are all the nodes of $R$.
Further, $c(\Delta)=c_0 - c(Q)$ for all $(0,q+1)$-paths $\Delta$, all $q$ of whose internal nodes are in $R$.
Therefore, condition (\ref{L2iii}) of Lemma \ref{L2} holds and Corollary \ref{C5} implies that $c_{ij}=0$ for all $(i,j) \in A(R \cup \{0\})$ and
$c_{i,q+1}=c_0$ for all $i \in R$. Replacing node $q+1$ by any other node in $V \setminus S$ (in the above argumentation), we obtain $c_{ij}=c_0$ for all 
$(i,j) \in (R:V \setminus S)$.

Next, consider for any arc $(i,j) \in A(V \setminus S)$ a tight \path $P$ that uses the arcs $(0,1),(1,2),(2,j)$ and skips node $i$.
Then the \path $P':= (P \setminus \{(0,1),(1,2),(2,j)\}) \cup \{(0,2),(2,i),(i,j)\}$ is also tight. Thus, we derive that $c_{ij}=0$ for all $(i,j) \in A(V \setminus S)$.
Further, from the tight \paths that starts with the arc $(0,1)$ and use some arc $(i,n)$ with $i \in V \setminus S$ we deduce $c_{in}=0$ for all those arcs $(i,j)$.
Moreover, from the tight \paths starting with the arc $(0,2)$ and ending with the arcs $(i,1),(1,n)$ for some $i \in V \setminus S$ we obtain $c_{i1}=0$ for 
$i \in V \setminus S$.
It is now easy to see that $c_{0i}=c_0$ for all $i \in V \setminus S$, $c_{jn}=0$ for all $j \in R$, and $c_{ij}=0$ for all $(i,j) \in (V \setminus S:R)$
(distinguish the cases $p=4$ and $p \geq 5$).
Therefore $\textbf{cx}=c_0$ is simply $c_0x((S:V \setminus S))=c_0$.
\end{proof}

\begin{Theorem}[cf. Theorem 13 of Hartmann and \oez \cite{HO}] \label{T13}
Let $S \subset V$ and $0,n \in S$. The one-sided min-cut inequality
\begin{equation} \label{osmincut1}
x((S: V \setminus S)) \geq x(\delta^+(l))
\end{equation}
is valid for $\Path$ for all $l \in V \setminus S$, and facet defining
for $\Path$ if and only if $|S| \geq p+1$ and $|V \setminus S| \geq 2$. 
\end{Theorem}

\begin{proof}
The one-sided min-cut inequality (\ref{osmincut1}) is valid, because all \paths that visits some node $l \in V \setminus S$ use at least one arc in $(S: V \setminus S)$.
If $|V \setminus S|=1$, then (\ref{osmincut1}) is the flow constraint $x(\delta^-(l))-x(\delta^+(l))=0$.
If indeed $|V \setminus S| \geq 2$ but $|S| \leq p$, then (\ref{osmincut1}) can be obtained by summing the min-cut inequality (\ref{mincut}) and the degree constraint $-x(\delta^+(l)) \geq -1$.

So suppose that $|S| \geq p+1$ and $|V \setminus S| \geq 2$. Let w.l.o.g. $l=1$ and set $R:=S \cup \{1\}$. By adding to (\ref{osmincut1}) the flow constraint $x(\delta^+(1))-x(\delta^-(1))=0$, it can be easily seen that (\ref{osmincut1}) is equivalent to 
\begin{equation} \label{mosmincut}
x((S: V \setminus R))- \sum_{i \in V \setminus R} x_{i1} \geq 0.
\end{equation}
Suppose that $\textbf{cx}=c_0$ is satisfied by every $\textbf{x} \in \Path$ that satisfies (\ref{mosmincut}) with equality.
By Corollary \ref{C4}, we may assume that $c_{in}=0$ for all $i \in V \setminus R$ and $c_{ij}=0$ for all arcs $(i,j)$ in some unbalanced 1-tree on $R$.
Condition (\ref{L2ii}) of Lemma \ref{L2} is satisfied; hence, from Corollary \ref{C5} follows that $c_{ij}=0$ for all $(i,j) \in A(R)$ which also implies that $c_0=0$.

Any \path that contains the arcs $(1,i),(i,n)$ for some $i \in V \setminus R$ and whose remaining arcs are in $A(R)$ satisfies (\ref{mosmincut}) with equality. 
Since $c_{in}=0$ and $c_a=0$ for all $a \in A(R)$, it follows that $c_{1i}=0$ for all $i \in V \setminus R$.
Now considering tight \paths that contain the arcs $(1,i),(i,j),(j,n)$ for some $(i,j) \in A(V \setminus R)$ and whose remaining arcs are in $A(R)$, 
we see that $c_{ij}=0$ for all $(i,j) \in A(V \setminus R)$. Further, the \paths that use the arcs $(1,i),(i,j)$ for $i \in V \setminus R, j \in S \setminus \{n\}$ 
and whose remaining arcs are in $A(R)$ yield $c_{ij}=0$ for all $(i,j) \in (V \setminus R: S \setminus \{n\})$.
Finally, considering for each $(i,j) \in (S:V \setminus R)$ and $k \in V \setminus R$ a tight \path that contains the arcs $(i,j),(j,1)$ 
and a tight \path that contains the arcs $(i,j), (j,k),(k,1)$, we see that $c_{j1}=c_{k1}$ for all $j,k \in V \setminus R$, 
$c_{ij}=c_{kl}$ for all $(i,j), (k,l) \in (S: V \setminus R)$, and
$c_{ij}+c_{k1}=0$ for all $(i,j) \in (S:V \setminus R)$, $k \in V \setminus R$.
Thus $\textbf{cx}=c_0$ is simply $c_{jk}x((S: V \setminus R))-
c_{jk}\sum_{i \in V \setminus R} x_{i1}= 0$ for some $(j,k) \in (S:V
\setminus R)$. 
\end{proof}

\begin{Theorem}[cf. Theorem 15 of Hartmann and \oez \cite{HO}] \label{T15}
Let $\langle R, S, T \rangle$ be a partition of $V$ and let $0,n \in
S$. The generalized max-cut inequality 
\begin{equation} \label{Roddmaxcut}
x((S:T)) + \sum_{i \in R} x(\delta^+(i)) \leq \lfloor (p+|R|)/2 \rfloor
\end{equation}
is valid for the $(0,n)-p$-path polytope $\Path$ for $p \geq 4$ and facet defining for $\Path$ if and only if
$p+|R|$ is odd, $|S \setminus \{n\}| > (p-|R|)/2$,
$|T| > (p-|R|)/2$, and
\begin{itemize}
\item[(i)] $p=|R|+3$, $|R| \geq 2$, and $|S|=3$, or
\item[(ii)] $p \geq |R|+5$.
\end{itemize}
\end{Theorem}

\begin{proof} 
\emph{Necessity}.
From $x(A)=p$ and $x((S:T)) \leq x((T:S))+x((T:R))$ we
  derive the inequality $2x((S:T)) + \sum_{i \in R} x(\delta^+(i))
  \leq p$. Adding the inequality $\sum_{i \in R} x(\delta^+(i)) \leq
  |R|$, dividing by two, and rounding down, we obtain (\ref{Roddmaxcut}).
When $p+|R|$ is even, then (\ref{Roddmaxcut}) is obtained
with no rounding, and hence it is not facet defining.
When $|S \setminus \{n\}| \leq (p-|R|)/2$ or $|T| \leq (p-|R|)/2$, then (\ref{Roddmaxcut}) is implied by degree constraints $x(\delta^(i)) \leq 1$.

Let $P$ be any \path and denote by $r$ the number of nodes in $R$ visited by $P$. Then $|v(P) \cap (S \setminus \{n\} \cup T)|=p-r$ and hence $\chi^P((S:T)) \leq (p-r)/2$. This in turn implies that there is no tight \path if $r \leq |R|-2$, where $|R| \geq 2$.
Now, when $p = |R|+3$ and $|S| \geq 4$, (\ref{Roddmaxcut}) is
dominated by nonnegativity constraints $x_{ij} \geq 0$ for $(i,j) \in
A(S \setminus \{0,n\})$.
Further, when $p=|R|+3$, $|S|=3$, and $|R|=1$, (\ref{Roddmaxcut}) is dominated
by the inequality (\ref{extra}).
Finally, when $p \leq |R|+1$, (\ref{Roddmaxcut}) is dominated by some
nonnegativity constraints, for example, $c_{in}=0$ for some $i \in T$.

\emph{Suffiency}. First we will show that (\ref{Roddmaxcut}) is facet
defining if $R= \emptyset$. In this the case, the resulting inequality
\begin{equation} \label{oddmaxcut}
x((S:T)) \leq \lfloor p/2\rfloor=q
\end{equation}
where $p=2q+1$, is called \emph{max-cut inequality}. First, we
show that (\ref{oddmaxcut}) is facet defining for $\Path$.
If $p=5$ and $|S \setminus \{n\}|=3$ or $|T|=3$, we will show that (\ref{oddmaxcut}) defines a facet using Theorem \ref{T8}.
The only primitive inequalities are those with $n=6$ and by
application of a convex hull code, we see that in this case
(\ref{oddmaxcut}) is facet defining for $\Path$. Moreover, 
(\ref{oddmaxcut}) is regular, since for each inner node $k$ there is a
non-tight \path that does not visit $k$. 
Without loss of generality, let $T= \{1,2,\dots,t\}$ and $S=
\{t+1,\dots,n,0\}$ for some $4 \leq t \leq n-4$. 

Suppose that $\textbf{cx}=c_0$ holds for all $\textbf{x} \in \Path$ satisfying (\ref{oddmaxcut}) with equality.
By Corollary \ref{C4}, we may assume that $c_{02}=1$, $c_{t+1,n}=0$, $c_{j1}=1$ for all $j \in S \setminus \{n\}$, and $c_{1i}=0$ for all $i \in T$.

First, consider any $(0,n)-2q$-path $P$ that alternates between nodes in $S$ and nodes in $T$, but does not visit node 1. Replacing
any arc $(i,j) \in P$ with $i \in S$, $j \in T$ by the arcs $(i,1),(1,j)$ we obtain a tight \path, and therefore
$c(P)-c_{ij}=c_0-1$ holds for all $(i,j) \in P \cap (S: T)$. This in turn
implies that $c_{ij}=1$ for all $(i,j) \in (S:T)$, since we have
$3 \leq t \leq n-3$ and $c_{02}=1$. Next, consider any tight \path
that uses arcs $(i,k), (k,j)$ for $i,j \in S \setminus \{0,n\}$, $k
\in T$ but does not visit node $l \in T$. Replacing node $k$ by node
$l$ yields another tight path which implies immediately
$c_{ik}+c_{kj}=c_{il}+c_{lj}$. Similarly we obtain
$c_{ki}+c_{il}=c_{kj}+c_{jl}$ and thus $c_{ik}+c_{ki}=c_{jl}+c_{lj}$
for all $i,j \in S \setminus \{0,n\}$ and $k,l \in T$. Since $t \geq
3$ and $c_{ik}=c_{jl}=1$, we see that there is some $\sigma$ with
$c_{ki}=\sigma$ for all $k \in T$, $i \in S \setminus \{0,n\}$. Now
consider any tight path that contains the arcs $(1,t+1),(t+1,n)$ and
does not visit some node $l \in T$. Replacing node $t+1$ by node $l$
yields another tight \path and hence
$c_{1,t+1}+c_{t+1,n}=c_{1l}+c_{ln}$. Since $c_{1,t+1}= \sigma$ and
$c_{t+1,n}=c_{1l}=0$, this implies $c_{ln}=\sigma$ for all $l \in T$,
$l \neq 1$. Of course, it follows also that $c_{1n}= \sigma$. 

Finally, any tight \path contains exactly one arc $(i,j) \in A(S) \cup
A(T)$, so $c_{ij}=c_0-q(1+ \sigma)$ for all $(i,j) \in A(S) \cup
A(T)$. Due to $c_{t+1,n}=0$, this implies that $c_{ij}=0$ for all
$(i,j) \in A(S) \cup A(T)$. Adding $\sigma$ times the equation
$x((S:T))-x((T:S))=0$, we see that $\textbf{cx}=c_0$ is equivalent
$x((S:T))=q$. This proves that (ref{oddmaxcut}) is also facet defining
when $0,n \in T$. 

When $R \neq \emptyset$, we prove the claim by showing that the
conditions of Theorem \ref{T9} hold for (\ref{oddmaxcut}). Since
$w=p-|R|$ is odd and $w \geq 5$, $x(S:T) \leq \lfloor w/2 \rfloor$
induces a facet of the $(0,n)-w$-path polytope defined on the digraph
$D^+=(V \setminus R, A(V \setminus R))$. Let us denote this inequality
by $\textbf{ax} \leq a_0$. It is easy to see that the auxiliary graph
$G_Z$ for the arc set $Z=\{(i,j)|a_{ij}=0\}$ is connected
(cf. \cite{HO}). Further, each tight $(0,n)-w$-path contains two arcs
$(i,j)$ and $(k,l)$ which are not adjacent and hence all arcs in $Z$
are in the same equivalency class with respect to $\mathbf{ax} \leq
a_0$. Since there are tight \paths with respect to (\ref{Roddmaxcut}) that
visit $|R|-1$ of the nodes in $R$, Theorem \ref{T9} implies that
(\ref{Roddmaxcut}) induces a facet of $\Path$ unless ($p=|R|+3, |R|
\geq 2, and |S|=3$).

Finally, suppose that $p=|R|+3$,$|R| \geq 2$, and $|S|=3$.
Without loss of generality, we may
assume that $S=\{0,1,n\}$, $2,3 \in R$, and $4,5 \in T$. 
Suppose that $\mathbf{cx}=c_0$ is satisfied by every $\mathbf{x}
\in \Path$ that satisfies (\ref{Roddmaxcut}) with equality.
By Corollary \ref{C4}, we may assume that
$c_{2j}=1$ for all $j \in  R$, $c_{i2}=0$ for all
$i \in T$,  $c_{32}=1$, $c_{21}=1$, $c_{1n}=0$, and $c_{04}=1$.
There are tight \paths that visits a node $l \in T$ followed by all
$|R|$ (or any $|r|-1$) nodes in $R$ and a node 1. Applying
Lemma \ref{L2}, we see that
\[
\begin{array}{rcll}
c_{lj} & = & \lambda + \pi_l - \pi_j & (j \in R) \\
c_{ij} & = & \lambda + \pi_i - \pi_j & (i,j \in R) \\
c_{im} & = & \lambda + \pi_i - \pi_m & (i \in R)
\end{array}
\]
for some $\lambda$, $\{\pi_j|j \in R\}$, $\pi_l$, and $\pi_1$. Let
w.l.o.g. $\pi_2=0$. Theorem \ref{T3} then implies that $\lambda =1$
and $\pi_j=0$ for all $j \in R$, $c_{i2}=0$ implies that $\pi_l=-1$,
and $c_{21}=1$ implies that $\pi_1=0$. Thus, $c_{ij}=1$ for all $(i,j) \in
A(R)$,  $c_{ij}=0$ for all $i \in T, j \in R$, and $c_{i1}=1$ for all
$i \in R$. Next, considering any (tight) \path $P$ that uses the arcs $(0,4),
(2,1), (1,n)$ and
visits all $|R|$ nodes in $R$ yields $c_0=|R|+1$. Replacing node 4 by
another node $j \in T$ yields $c_{0j}=1$ for all $j \in
T$. Next, consider any tight \path $P$ that uses the arcs $(0,i),
(i,j), (j,1)$ for some $i,j \in R$. Then the \path $P':=(P \setminus
\{(0,i),(i,j),(j,1)\}) \cup \{(0,j),(j,i),(i,1)\}$ is also tight, and
hence, $c_{0i}=c_{0j}$ for all $i,j \in R$. Denote this common value
by $\sigma$. From the tight \paths that visits the nodes 1 and $t$ for
some $t \in T$ and all nodes in $R$, we derive $c_{ij}=1- \sigma$ for
all $i \in R,j \in T$. Now it is easy to see that $c_{in}=1+ \sigma$
for all $i \in T$. Considering any tight \path that uses the arcs
$(0,2),(2,1),(1,4), (4,3)$, and $(m,n)$ for an appropriate $m \in R$
yields $\sigma = 0$. Thus, $c_{0i}=0$ and $c_{in}=1$ for all $i \in
R$, $c_{1j}=1$ for all $j \in T$, and $c_{ij}=1$ for all $i \in R, j
\in T$. Determining the coefficients of the remaining arcs is an easy
task. So we see that $\mathbf{cx}=c_0$ is simply (\ref{Roddmaxcut}).
\end{proof}

\begin{Theorem} \label{T16}
Let $\langle R, S, T \rangle$ be a partition of $V$ and let $0,n \in
T$. The generalized max-cut inequality 
\begin{equation} \label{Roddmaxcut1}
x((S:T)) + \sum_{i \in R} x(\delta^+(i)) \leq \lfloor (p+|R|)/2 \rfloor
\end{equation}
is valid for the $(0,n)-p$-path polytope $\Path$ for $p \geq 4$ and
facet defining for $\Path$ if and only if $p+|R|$ is odd, $|S| > (p-|R|)/2$,
$|T \setminus {0}| > (p-|R|)/2$, and 
\begin{itemize}
\item[(i)] $p=|R|+3$, $|R| \geq 2$, and $|T|=3$, or
\item[(ii)] $p \geq |R|+5$.
\end{itemize} \hfill $\Box$
\end{Theorem}

\begin{Theorem} \label{Tevencut1}
Let $\langle R, S, T \rangle$ be a partition of $V$, let $0 \in S$, and let $n \in T$. The generalized max-cut inequality
\begin{equation} \label{Revenmaxcut}
x((S:T)) + \sum_{i \in R} x(\delta^+(i)) \leq \lfloor (p+|R|+1)/2 \rfloor
\end{equation}
is valid for the $(0,n)-p$-path polytope $\Path$ for $p \geq 4$ and facet defining for $\Path$ if and only if
$p+|R|$ is even, $p \geq |R|+4$, $|S| > (p-|R|)/2$, and $|T| > (p-|R|)/2$.
\end{Theorem}

\begin{proof}
From the equation $x(A)=p$ and the inequality $x((S:T)) \leq x((T:S))+x((T:R))+1$ we derive the inequality $2x((S:T))+ \sum_{i \in R}x(\delta^+(i)) \leq p+1$. Adding the inequality $\sum_{i \in R} x(\delta^+(i)) \leq |R|$, dividing by two, and rounding down yields (\ref{Revenmaxcut}). If $p+|R|$ is odd
we obtain (\ref{Revenmaxcut}) without rounding and hence it is not facet defining for $\Path$.
When $|S| \leq (p-|R|)/2$ or $|T| \leq (p-|R|)/2$, (\ref{Revenmaxcut})
is dominated by degree constraints $x(\delta^+(j)) \leq 1$. 
Furthermore, we have to show that (\ref{Revenmaxcut}) is not facet
defining if $p \leq |R|+2$. 
When $R= \emptyset$, it is clear. Otherwise consider any \path $P$ and denote the number of nodes in $R$ visited by $P$ by $r$. It is easy to see that
$P$ is tight only if $r \geq |R|-1$. For the sake of contradiction, assume that $p \leq |R|$ and $P$ is tight. Then we have $r=|R|-1$ and thus $p=|R|$
which implies $\lfloor (p+|R|+1)/2 \rfloor = |R|$. But $\chi^P((S:T))
+ \sum_{i \in R} chi^P(\delta^+(i))=|R|-1$, so $P$ is not tight, a contradiction.
Hence, the only possibility is that $p=|R|+2$. Now, $p=|R|+2$ implies that $|S|,|T| \geq 2$ and $\lfloor (p+|R|+1)/2 \rfloor = |R|+1$. But then
(\ref{Revenmaxcut}) is dominated by the nonnegativity constraints $x_{ij} \geq 0$ for all $(i,j) \in A(S) \cup A(T)$.

First, we show that (\ref{Revenmaxcut}) is facet defining when $R=
\emptyset$. In this case, $p$ is even and (\ref{Revenmaxcut}) is the
max-cut inequality
\begin{equation} \label{evenmaxcut}
x((S:V \setminus S)) \leq \lfloor (p+1)/2 \rfloor =p/2.
\end{equation}
If $p=4$ and $|S|=3$ or $|V \setminus S|=3$, we will show that
(\ref{evenmaxcut}) defines a facet of $\Path$ using Theorem
\ref{T8}. The only primitive members of family (\ref{evenmaxcut}) with
$p=4$ are those with $|S|=|V \setminus S|=3$. Inequality
(\ref{evenmaxcut}) is obviously regular and using a convex hull code,
we see that (\ref{evenmaxcut}) defines a facet of $\Path$. Moreover,
all \pbowties tied at an inner node satisfy (\ref{evenmaxcut}). 

If $p \geq 6$ suppose that the equation $\textbf{cx}=c_0$ is satisfied
by every $\textbf{x} \in \Path$ that satisfies ((\ref{evenmaxcut}))
with equality. Let w.l.o.g. $1,2 \in V \setminus S$. By Corollary
\ref{C4}, we may assume that $c_{02}=1$, $c_{i1}=1$ for all $i \in S$,
and $c_{1j}=0$ for all $j \in V \setminus S$, $j \neq 1$. Since $|S|,
|V \setminus S| \geq 4$, we can apply the same argumentation as in the
proof to Theorem \ref{T15}. Thus $c_{ij}=1$ for all $(i,j) \in (S: V
\setminus S)$, $c_{ij}=\sigma$ for all $(i,j) \in (V \setminus (S \cup
\{n\}):S \setminus \{0\})$, for some $\sigma$, and $c_{ij}=0$ for all
$(i,j) \in A(S) \cup A(T)$. Evaluating the cost of tight \paths yields
$c_0= \frac{p}{2}+ (\frac{p}{2}-1) \sigma$ which implies that
$\textbf{cx}=c_0$ is the equation $x((S:V \setminus S)) + \sigma x((V
\setminus S:S))=\frac{p}{2}+\sigma(\frac{p}{2}-1)$. Adding $\sigma$
times the equation $x((S:V \setminus S))-x((V \setminus S:S))=1$, we
see that (\ref{evenmaxcut}) is equivalent to $x((S:V \setminus
S))=p/2$. 

Applying Theorem \ref{T9} to the $(0,n)-w$-path polytope defined on
the digraph $D^*=(V \setminus R, A(V \setminus R))$, where $w=p-|R|$,
proves that (\ref{Revenmaxcut}) is facet defining for $\Path$ even for
$R \neq \emptyset$. 
\end{proof}

\begin{Theorem} \label{Tevencut2}
Let $\langle R, S, T \rangle$ be a partition of $V$, let $0 \in T$,
and let $n \in S$. The generalized max-cut inequality 
\begin{equation} \label{Revenmaxcut2}
x((S:T)) + \sum_{i \in R} x(\delta^+(i)) \leq \lfloor (p+|R|-1)/2 \rfloor
\end{equation}
is valid for the $(0,n)-p$-path polytope $\Path$ for $p \geq 4$ and
facet defining for $\Path$ if and only if $p+|R|$ is even, $p \geq
|R|+4$, $|S| > (p-|R|)/2$, and $|T| > (p-|R|)/2$.
\hfill $\Box$
\end{Theorem}

\begin{Remark}
If $R = \emptyset$, inequality \emph{(\ref{Revenmaxcut2})} is equivalent to
the inequality 
\[
x((T:S)) \leq  \lfloor (p+1)/2 \rfloor,
\]
since in this case holds the equation $x((S:T))=x((T:S))-1$.
\end{Remark}

\begin{Theorem}
Let $\emptyset \neq T = V \setminus \{0,1,2,3,n\}$. The
inequality
\begin{equation} \label{extra}
\begin{array}{rcl}
x_{03}-x_{3n}+3x_{12}-x_{21}+2x_{13}-2x_{31}-2x_{2n}+2x((T:\{3\})) & & \\
+x(A(T))+x((\{1\}:T))-x((T:\{1\}))+x((T:\{2\}))-x((\{2\}:T) & \geq & 0
\end{array}
\end{equation}
is facet defining for $P_{(s,t)- \scriptsize{path}}^4(D)$.
\end{Theorem}

\begin{proof}
When $|T|=1$, the claim can be verified with a convex hull code. For
$|T| \geq 2$ we apply Theorem \ref{T8}. 
\end{proof}

\subsection{Jump inequalities}

Dahl and Gouveia \cite{Dahl} introduced a class of valid inequalities
for the directed hop-constrained shortest path problem 
(the problem of finding a minimum $(0,n)$-path with at most $p$ arcs)
they called \emph{jump and lifted jump inequalities}. Given a partition 
 $\langle S_0,S_1,S_2,\dots,S_p,S_{p+1} \rangle$ of $V$ into $p+2$
 node sets, where $S_0=\{0\}$ and $S_{p+1}=\{n\}$, 
these inequalities encode the fact that a
 $(0,n)$-path $P$ of length at most $p$ must make at least one "jump" 
from a node set $S_i$ to a node set $S_j$, with $j-i \geq 2$.
Transferring them to the \path polytope and lifting them (see
\cite{Dahl})  we can give a sufficient condition for them to be facet
defining for $\Path$. But it seems to be hard to give a complete
classification of the jump inequalities.

\begin{Theorem}
Let $\langle S_0,S_1,S_2,\dots,S_p,S_{p+1} \rangle$ be a partition of
$V$, where $S_0=\{0\}$ and $S_{p+1}=\{n\}$. The jump inequality 
\begin{equation} \label{jump}
\sum_{i=0}^{p-1} \sum_{j=i+2}^{p+1} x((S_i:S_j)) - x((S_{p-1}
\cup S_p : S_1 \cup S_2)) \geq 1
\end{equation}
is facet defining for the \path polytope
$\Path$ if $|S_i| \geq 2$ for $i=1,\dots,p$.
\end{Theorem}

\begin{proof}
We refer to an arc $(i,j)$ as forward arc if $(i,j) \in (S_k:S_l)$ for
some $k<l$ and as backward arc if $(i,j) \in (S_q:S_r)$ for some
$q>r$. We say, the \path $P$ makes a ``jump'' with respect to
(\ref{jump})if $P$ uses an arc $(i,j) \in (S_k:S_l)$ for some $0 \leq
k < l \leq p+1$ with $l \geq k+2$. 

 The jump inequality (\ref{jump}) is valid for $\Path$, since it is
 valid for the path polytope $\Pathleq$ which is the convex hull of all
 incidence vectors of simple $(0,n)$-paths with at most $p$ arcs (see
 \cite{Dahl}). 

To show that (\ref{jump}) is facet defining for $\Path$, we apply
Theorem \ref{T8}. So we have to verify that the conditions of Theorem
\ref{T8} hold for (\ref{jump}), when $|S_i|=2$ for $i=1,\dots,p$, that
is, when $n=2p+1$.
In the sequel, let $\mathbf{dx} \geq 1$ be such an jump inequality.

Let $B=P \cup C$ be any \pbowtie, where $C$ is a simple cycle and $P$
is a simple $(0,n)$-path. Since $|P| \leq p$, $d(P) \geq 1$. When
$d(C) \geq 0$, it follows $d(B) \geq 1$, too. Otherwise $d(C)=-1$
and $C$ is a cycle in 
\[\left( \bigcup_{j=2}^{p-2} (S_j:S_{j+1})  \right) \cup (S_{p-1}:S_2),\]
since $|C| \leq p-2$. Thus, the cardinality of $C$ is equal to $p-2$
and $P$ is a $(0,n)-2$-path that makes two ``jumps''.
Therefore,  the jump inequality $\mathbf{dx} \geq 1$ is satisfied by
all \pbowties.

Further,  $\mathbf{dx} \geq 1$ is regular, since to each
internal node $k$ there exists a non-tight \path that does not visit
node $k$.

It remains to be shown that $\mathbf{dx} \geq 1$ is facet defining for
$\Path$. Without loss of generality, let $S_i=\{i,p+i\}$ for
$i=1,\dots,p$.
When $p=4$ or $p=5$, the inequality $\mathbf{dx} \geq 1$ can be seen
facet defining using a convex hull code. So let $p \geq 6$.
Suppose that $\mathbf{cx}=c_0$ is satisfied by every $\mathbf{x} \in \Path$
that satisfies (\ref{jump}) with equality. Denoting by $P$ the
$(0,2p+1)$-path $(0,\dots,p,2p+1)$, we may
assume by Corollary \ref{C5} that $c(P)=0$, $c_{0,p+1}=0$, and
$c_{i,p+i}=0$ for $i=1,\dots,p$. Substituting two connected arcs
$(i,j), (j,k) \in P$ by the arc $(i,k)$, we see that $c_{m-1,m+1}=c_0$
for $m=1,\dots,p-1$, and $c_{p-1,2p+1}=c_0$. Next, replacing three
connected arcs $(i,j),(j,k),(k,l) \in P$ with $i > 0$ by the arcs
$(i,p+i),(p+i,l)$, we see that $c_{2p-2,2p+1}=c_0$ and
$c_{p+i,i+3}=c_0$ for $i=1,\dots,p-3$.
Further, replacing in these \paths node $i$ by node $p+i-1$ (for $i
\geq 2$) yields $c_{m,m+1}=0$ for $m=p+1,\dots,2p-3$ and considering
successively the \paths 
\[
(0,p+1,4,\dots,q,p+q,\dots,2p+1)
\]
for $q=p,\dots,4$, we see that even $c_{m,m+1}=0$ for
$m=p+1,\dots,2p$, since $p \geq 6$.
We can now easily deduce that $c_{i,p+i+1}=c_{p+i,i+1}=c_{p+i,i}=0$
for $i=1,\dots, p$, $c_a=c_0$ for all $a \in (S_i:S_{i+2})$
($i=0,\dots,p-1$), and $c_a=c_0$ for all $a \in (S_i:S_{i+3})$
($i=0,\dots,p-2$). Furthermore, for each arc $a \in (S_i:S_{i+4})$,
$i=0,\dots,p-3$, there is a tight \path containing $a$ that does not
use any backward arc, which implies that $c_a=c_0$ for all those arcs $a$.
Moreover, for each arc $a \in (S_m:S_{m-1})$ there is a tight \path
that uses $a$, makes a jump from $S_i$ to $S_{i+4}$ for some $i$, and
does not use any further backward arcs. Hence, $c_a=0$ for all
$a \in (S_m:S_{m-1})$, $m=2,\dots,m$.
It is now easy to see that the remaining coefficients can be
determined as required, and therefore, $\mathbf{cx}=c_0$ is simply
$c_0 \mathbf{dx}=c_0$.
\end{proof}

\subsection{Cardinality-path inequalities}
The cardinality-path inequalities were originally formulated for the
cardinality constrained circuit polytope. They say that a (undirected)
simple cycle of cardinality at most $p$ never uses more edges of a
(undirected) simple path $P$ of cardinality $p$ than internal nodes of
$P$. This idea can be transferred to the \path polytope. Before
stating the next theorem we introduce two notations. For any simple
path $P$ we denote its internal nodes by $\dot{P}$. Furthermore, we define
$\bid(P):= P \cup \{(i,j)|(j,i) \in P\}$.

\begin{Theorem} \label{T20}
Let $s,t$ be internal nodes and $P$ be a $(s,t)$-path of length $p-1$.
The cardinality path inequality
\begin{equation} \label{cardpath}
\sum_{i \in \dot{P}} x(\delta^-(i)) - x(\bid(P)) \geq 0
\end{equation}
is valid for the \path polytope $\Path$ and induces a facet of $\Path$
if and only if $p \in \{4,5\}$ and $n \geq p+2$ or $p \geq 6$ and $n
\geq 2p -3$.  
\end{Theorem}

\begin{proof}
Without loss of generality, let $P=(1,2,\dots,p)$.

\emph{Necessity.} When $p \in \{4,5\}$ and $n=p+1$, (\ref{cardpath})
can be seen not to induce a facet using a convex hull code. When $p \geq 6$ and
$p+1 \leq n \leq 2p-4$, (\ref{cardpath}) is dominated by the
nonnegativity constraints $x_{2,p-1} \geq 0$ and $x_{p-1,2} \geq 0$. 

\emph{Suffiency.}
When the conditions in Theorem \ref{T20} are satisfied and the
cardinality of the node set $S:= \{1,p,p+1,\dots,n-1\}$ is at most 4,
(\ref{cardpath}) can be seen to induce a facet using a convex hull code. 
So suppose that $|S| \geq 5$ and $\mathbf{cx}=c_0$ is satisfied by
every $\mathbf{x} \in \Path$ that satisfies (\ref{cardpath}) with
equality. 
By Corollary \ref{C5} we may assume that $c_{j,j+1}=0$ for
$j=1,\dots,p-2$, $c_{0,n-1}=c_{n-1,n}=0$, and $c_{ij}=0$ for all arcs
$(i,j)$ in some unbalanced 1-tree on $S$.  

For any four distinct nodes $i \in S \cup \{0\}$, $j,k \in S$, and $l \in S \cup \{n\}$ there is a tight \path that uses the arcs $(i,k),(k,j)$ and skips node $l$. Replacing node $k$ by node $l$ yields another tight \path and thus $c_{ik}+c_{kj}=c_{il}+c_{lj}$. Using Corollary \ref{C5} we obtain 
$c_{ij}=0$ for all $(i,j) \in A(S \cup \{0,n\})$ and therefore also $c_0=0$. 

In the following we distinguish the three cases $p=4$, $p=5$, and $p \geq 6$.

\vspace{0.5cm}
\noindent CASE 1: $p=4$

From the \fourpaths $(0,5,1,2,n)$ and $(0,1,2,3,n)$ we derive $c_{2n}=c_{3n}=0$ and from the \fourpaths $(0,1,2,i,n)$ for
$i=p,\dots,n-1$ we derive $c_{2i}=0$.

Next, considering the \fourpaths $(0,5,4,3,n)$ and  $(0,4,3,2,n)$ yields $c_{43}=c_{32}=0$. Hence, we can also deduce that $c_{3j}=0$ for all $j \in S \setminus \{4\}$.

Further, from all tight $(0,n)-4$-paths that use the arc $(3,4)$ we deduce that $c_{ij}+c_{34}=0$ for $(i,j) \in \{(0,2),(0,3),(1,3)\} \cup (S: \{3\})$. It follows analogously that $c_{kl}+c_{21}=0$ for all $(k,l) \in \{(0,2),(0,3),(4,2)\} \cup (S: \{2\})$.
In particular, $c_{02}+c_{21}=c_{02}+c_{34}=0$ which implies that $c_{21}=c_{34}$ and hence, $c_{ij}+c_{kl}=0$ for all
$(i,j) \in \{(1,3),(4,2)\} \cup (S \cup \{0\}: \{2,3\})$ and $(k.l) \in \{(2,1),(3,4)\}$. So $\mathbf{cx}=c_0$ is obviously equivalent to (\ref{cardpath}).

\vspace{0.5cm}
\noindent CASE 2: $p=5$

This case can be carry out similar as the case $p=4$; so we omit this part of the proof.

\vspace{0.5cm}
\noindent CASE 3: $p \geq 6$

From the \path $(0,\dots,p-1,n)$ we derive that $c_{p-1,n}=0$. Further, setting $T:= \{3,\dots,p-2\}$, it can be easily seen that $c_{ij}=0$ for all $i \in T, j \in (S \setminus \{1\}) \cup \{n\}$. 
Next, for any arc $(i,j) \in (\dot{P} \setminus \{p-1\}: S \cup \{n\} \cup \{(p-1,n)\})$ there is a tight \path that uses the arcs $(i,j)$ and $(k,k+1)$ for $k=1,\dots,i-1$ and whose remaining arcs are in $A(S \cup \{0,n\})$. Hence, $c_{ij}=0$ for all those arcs $(i,j)$.
Further, from the \path $(0,\dots,p-3,p,p-1,n)$ we derive that $c_{p,p-1}=0$. Moreover,
for any node $i \in S \setminus \{1\}$ there is a tight \path that
uses the arcs $(0,1),(1,2),(2,i),(p,p-1),(p-1,n)$ and whose remaining
arcs are in $A(S)$. Thus, $c_{2i}=0$ for all $i \in S \setminus
\{1\}$. Considering further tight \paths on node set $S \cup
\{0,2,p-1,n\}$, we see that also $c_{p-1,i}=0$ for all $i \in S
\setminus \{p\}$ and $c_{2n}=0$. Finally, considering successively the
\paths $(0,\dots,i-2,p,p-1,\dots,i,n)$
for $i=p-2,\dots,2$, we find that $c_{i+1,i}=0$ for $i=2,\dots,p-2$.

It remains to be shown that $c_{21}=c_{p-1,p}=\sigma$ and
$c_{ij}=-\sigma$ for all arcs $(i,j)$ in $\bigcup_{k=2}^{p-1}
\delta^-(k) \setminus \bid(P) $ for some
$\sigma$. From the two tight \paths $(0,4,5,\dots,p+2,n)$ and
$(0,4,3,2,1,p+1,p+2,\dots,n)$ we derive that
$c_{21}=c_{p-1,p}$. Denote this common value by $\sigma$.
Since to each arc $(i,j) \in \bigcup_{k=2}^{p-1}
\delta^-(k) \setminus \bid(P)$ there is a tight \path that uses either
the arc $(2,1)$ or $(p-1,p)$ and therefore, $c_{ij}=-\sigma$ for
all those arcs $(i,j)$. Thus, $\mathbf{cx}=c_0$ is simply
\[\sigma x(bid(P))- \sigma \sum_{i \in V(\dot{P})} x(\delta^-(i))=0.\]
\end{proof}

\section{Facets of related polytopes}

In this section, we derive facet defining inequalities for related
polytopes from facet defining inequalities for the \path polytope. We
exploit three tools to do this; the first is Theorem \ref{T0} which
can be applied  to derive facets for the $p$-cycle polytope. The two
other tools were already mentioned in Hartmann and \oez \cite{HO}.
They showed that the undirected counterpart
$\mathbf{\bar{c}y} \leq c_0$  of a
symmetric inequality  $\mathbf{cx} \leq c_0$  is facet inducing for the (undirected)
$p$-circuit polytope $P_C^p(K_n)$ if  $\mathbf{cx} \leq c_0$ is
facet inducing for $P_C^p(D_n)$. Here,  $\mathbf{cx} \leq c_0$ is
called \emph{symmetric} if $c_{ij}=c_{ji}$ for all $i < j$ and
the induced inequality $\mathbf{\bar{c}y} \leq c_0$ for $P_C^p(K_n)$
is defined by $\bar{c}_{ij}=c_{ij}=c_{ji}$ for all $i < j$.
This concept can be adapted to the directed and undirected path polytopes in
a modified version. We refer to \ref{up}. The third tool can be applied to
the undirected/directed \path or $p$-cycle polytopes (basic polytopes),
when relaxing the cardinality constraint $x(B)=p$ to $x(B) \geq p$ or
$x(B) \leq p$, where $B$ is the ground set (the arc set or edge set). 
The resulting upper and lower polytopes have one dimension more than their basis polytopes, respectively, and this fact can be exploited to lift facets of the
basis polytope into facets of the related upper and lower polytopes (see \ref{relax}).

We illustrate the three tools by examples in the next subsections. In
\ref{up}, we apply not only the second tool, but also give a short
polyhedral analysis of the undirected counterpart of the \path
polytope.

\subsection{New facets of the directed $p$-cycle polytope}
Applying Theorem \ref{T0} to Theorems \ref{Tevencut1} and
\ref{Tevencut2} we obtain some new facet defining
inequalities for the directed $p$-cycle polytope $P_C^p(D_n)$.

\begin{Corollary}
Let $\langle \{j\}, R, S, T \rangle$ be a partition of $V$.
 The inequality
\begin{equation} \label{Revenmaxcut1}
x((S:\{j\}))+x((\{j\}:T))+x((S:T)) + \sum_{i \in R} x(\delta^+(i))
 \leq \lfloor (p+|R|+1)/2 \rfloor
\end{equation}
defines a facet of the $p$-cycle polytope $P_C^p(D_n)$ if
$p+|R|$ is even, $p \geq |R|+4$, $|S| > (p-|R|)/2-1$, and $|T| >
(p-|R|)/2-1$.
\hfill $\Box$
\end{Corollary}

\begin{Corollary}
Let $\langle \{j\}, R, S, T \rangle$ be a partition of $V$.
 The inequality
\begin{equation} \label{Revenmaxcut3}
x(\delta^+(r))+x((S:T)) + \sum_{i \in R} x(\delta^+(i))
 \leq \lfloor (p+|R|+1)/2 \rfloor
\end{equation}
defines a facet of the $p$-cycle polytope $P_C^p(D_n)$ if
$p+|R|$ is even, $p \geq |R|+4$, $|S| > (p-|R|)/2-1$, and $|T| >
(p-|R|)/2-1$.
\hfill $\Box$
\end{Corollary}

\subsection{Facets of the undirected $(0,n)-p$-path polytope} \label{up}

The undirected  $(0,n)-p$-path polytope $\UPath$ is the symmetric counterpart
of the directed \path polytope $\Path$. Here, $K_{n+1}=(V,E)$ denotes the
complete graph on node set $V=\{0,\dots,n\}$.
Table 2 gives linear descriptions of $\UPathone$ and
$\UPathtwo$. The complete polyhedral analysis of the \upath polytope
$\UPaththree$ begins with the next theorem and afterwards we will turn
to the \upath polytopes $\UPath$ with $4 \leq p \leq n-1$.

\begin{Theorem}
Let $K_{n+1}=(V,E)$ be the complete graph on node set
$V=\{0,\dots,n\}$. Then
\[\dim \UPaththree = |E|-n-2.\]
\end{Theorem} 

\begin{proof}
First note that each internal edge $e=[i,j]$ corresponds to two
incidence vectors $P^{(i,j)}$ and $P^{j,i}$ of $[0,n]-3$-paths as
follows: 
$P^{(i,j)}= \chi^{[0,i],[i,j],[j,n]}$ and 
$P^{(j,i)}= \chi^{[0,j],[j,i],[i,n]}$.
Consider the points $P^{(k,n-1)}$, $P^{(n-1,k)}$ for
$k=1,\dots,n-2$ and $P^{(i,j)}$ for $1 \leq i < j \leq n-2$.
It is easy to see that these $|E|-n-1$ points are linearly independent
and thus, $\dim \UPaththree \geq |E|-n-2$.

Next, all incidence vectors of $[0,n]-3$-paths satisfy the following system
of linearly independent equations:
\begin{eqnarray}
y_{0n} & = & 0, \label{101} \\
y(\delta(0)) & = & 1, \label{102}  \\
y(\delta(n)) & = & 1, \label{103}  \\
y(\delta(i)) -2(y_{0i}+y_{in}) & = & 0, \hspace{1cm} i=1,\dots,n-1, \label{104}
\end{eqnarray}
where $\delta(j)$ denotes the set of edges which are incident with
node $j$ and $y(F)= \sum_{e \in F} y_e$ for any $F \subseteq E$.
This implies that $\dim \UPaththree \leq |E|-n-2$, which completes the proof.
\end{proof}

\begin{Remark}
Adding the
equations (\ref{102})-(\ref{104}), subtracting two times \emph{(\ref{101})}, and dividing by two, yields the
equation
\begin{equation} \label{106}
\sum_{i=1}^{n-2} \sum_{j=i+1}^{n-1} y_{ij}=1.
\end{equation}
\end{Remark}

In the next theorem, $\delta_{in} \in \{0,1\}$ for $i=1,\dots,n-1$.

\begin{Theorem} \label{T100}
A complete and nonredundant linear description of the $[0,n]-3$-path
polytope $\UPaththree$ 
is given by the equations \emph{(\ref{101})-(\ref{104})}, the nonnegativity
constraints $y_{ij} \geq 0$ for $1 \leq i < j \leq n$, and the inequalities
\begin{equation} \label{105}
\sum_{i=1}^{n-1} \delta_{in} y_{in} + \sum_{i=1}^{n-2}
\sum_{j=i+1}^{n-1} \left \lfloor \frac{2- \delta_{in}- \delta_{jn}}{2}
\right \rfloor y_{ij} \leq 1
\end{equation}
for all $(n-1)$-tupels $(\delta_{1n},\dots, \delta_{n-1,n})$
satisfying $1 \leq \sum_{i=1}^{n-1} \delta_{in} \leq n-2$.
\end{Theorem}

\begin{table}
\noindent Table 2. Polyhedral analysis of $\UPathone$ and $\UPathtwo$.

\vspace{0.2 cm}
\noindent
\begin{tabular}{|c|c|rcll|} \hline
& & & & &\\
\rb{$p$} & \rb{Dimension} & \multicolumn{4}{c|}{\rb{Complete linear description}} \\ \hline 
      &  & $y_{0n}$ & $=$ & $1$ & \\
\rb{$1$} & \rb{0} & $y_{ij}$ & $=$ & $0$ \hspace{1.6 cm} & $\forall \;  [i,j] \in E \setminus \{[0,n]\}$ \\ \hline \hline 
& & $y_{0n}$ & $=$ & $0$ & \\ 
& & $y(\delta(0))$ & $=$ & $1$ & \\
$2$ & $n-2$ & $y_{0i} - y_{in}$ & $=$ & $0$ & $i=1,\dots,n-1$ \\
& & $y_{0i}$ & $\geq $ & $0$ & $i=1,\dots,n-1$ \\
& & $y_{ij}$ & $=$ & $0$ & $1 \leq i < j \leq n-1$ \\ \hline
\end{tabular}
\end{table}

\begin{proof}
\emph{Validity.} 
Let $\mathbf{cy} \leq 1$ be some inequality of family (\ref{105}).
The edge set of the support graph $G=(V,F)$, defined by
$F:=\{e \in E| c_e=1\}$, decomposes into two disconnected subsets
$F^n:=\{[i,n] \in F| \delta_{in}=1\}$ and $F^{\neg n}:=F \setminus
F^n$, and as is easily seen, each $[0,n]-3$-path $P$ uses
at most one edge of $F$ in the subgraph $G \subset K_{n+1}$. Hence,
$\mathbf{cy} \leq 1$ is valid for $\UPath$.

\emph{Nonredundancy.} 
Since the equations (\ref{101})-(\ref{104}) are linearly independent,
they induce a nonredundant description of the lineality space of
$\UPaththree$.

Next, we prove that the inequalities given in Theorem \ref{T100}
are nonredundant by showing that the set of induced faces is an anti-chain. Let $F_1$ and $F_2$ be from two different
inequalities induced faces of
$\UPaththree$. When $F_1$ and $F_2$ are induced by nonnegativity
constraints, they are clearly not contained into each other. If only
one of them is induced by a nonnegativity constraint $y_{ij} \geq 0$
($1 \leq i < j \leq n$), say $F_1$,
it follows immediately that $F_2 \not\subset F_1$. 
Since $|V(F^{\neg n})| \geq 2$, there is also a point $P^{(k,l)}$ in
$F_1$ that is not in $F_2$ and thus,  $F_1 \not\subset F_2$. 

Finally, let both faces not induced by nonnegativity
constraints. Denote the edge sets of the support graphs corresponding to
$F_1$ and $F_2$ by $E_1$ and $E_2$, respectively. Since $E_1
\not\subset E_2$ and $E_2 \not\subset E_1$, it follows also that 
$F_1 \not\subset F_2$ and $F_2 \not\subset F_1$.

\emph{Completeness.}
We will show that each facet defining inequality $\mathbf{cx} \leq c_0$ for $\UPaththree$ is
equivalent to a nonnegativity constraint $y_{ij} \geq 0$ ($1 \leq i <
j \leq n$) or an inequality of family (\ref{105}).

Adding appropriate multiples of the equations (\ref{101})-(\ref{104}),
we see that $\mathbf{cy} \leq c_0$ is equivalent to an inequality 
$\mathbf{dy} \leq d_0$ with
\begin{itemize}
\item[(i)] $d_{0i}=0$ for $i=1,\dots,n$, 
\item[(ii)] $d_{zn}=0$ for some internal node $z$,
\item[(iii)] $d_{uw}=0$ for some internal edge $[u,w]$, and
\item[(iv)] $d_{ij} \geq 0$ for $1 \leq i < j \leq n$. 
\end{itemize}
This immediately implies that $d_0 > 0$ and $0 \leq d_e \leq d_0$ for
all $e \in E$. 

Next, we will show that $d_e \in \{0,d_0\}$ for all $e \in E$.
Suppose, for the sake of contradiction, that $M:= \{[i,j] \in E|
0 < d_{ij} < d_0\} \neq \emptyset$. Assuming that there is some
internal edge $[k,l] \in M$ with $[k,n], [l,n] \notin M$, we see that
$d_{kn}=d_{ln}=0$, since $d_{kl}+d_{ln} \leq d_0$ and $d_{kl}+d_{kn}
\leq d_0$. Thus,  $\mathbf{dy} \leq d_0$ is
dominated by the inequality $\mathbf{\tilde{d}y} \leq d_0$, where
$\tilde{d}_{kl}=d_0$ and $\tilde{d}_e=d_e$ for all $e \in E \setminus
\{[k,l]\}$. Assuming that there is some edge $[m,n]$ such that $[i,m]
\notin M$ for all internal nodes $i \neq m$,
yields $d_{im}=0$ for all internal nodes $i \neq m$. Therefore
$\mathbf{dy} \leq d_0$ is dominated by the inequality $\mathbf{d'}y
\leq d_0$, where $d'_{mn}=d_0$ and $d'_e=d_e$ for all $e \in E
\setminus \{[m,n]\}$. So we may assume in the sequel: 
\begin{itemize}
\item[(a)]  $[i,n] \in M$ or $[j,n] \in M$ for each internal edge
  $[i,j] \in M$; 
\item[(b)]  for each edge $[k,n] \in M$ there is an internal edge
  $[i,k] \in M$. 
\end{itemize} 
In particular, we deduce that $M \cap \{[i,j]| 1 \leq i < j \leq n-1\}
\neq \emptyset$ and $M \cap \{[i,n] | 1 \leq i \leq n-1\} \neq \emptyset$.

Let $d_{rs}$ be the minimum over all
edges in $M \cap \{[i,j] | 1 \leq i < j \leq n-1\}$ and $d_{vn}$ be
the minimum over all edges in $M \cap \{[i,n]| 1 \leq i \leq n-1\}$.
We now construct two different inequalities $\mathbf{ay} \leq a_0$ and
$\mathbf{by} \leq b_0$ that together imply $\mathbf{dy} \leq d_0$.
The coefficients of the both inequalities we set as follows:
\[
\begin{array}{ccccll}
a_0 & = & b_0 & = & d_0, \\
a_{ij} & = & b_{ij} & = & d_{ij} & \forall \; [i,j] \in E \setminus M,\\
a_{ij} & = &        &   & d_{ij}-d_{rs} & \mbox{for } 1 \leq i < j
\leq n-1,\\
a_{kn} & = &        &   & d_{kn} + d_{rs} & \mbox{for } 1 \leq k \leq n-1,\\
       &   & b_{ij} & = & d_{ij} + d_{vn} & \mbox{for } 1 \leq i < j
\leq n-1,\\
       &   & b_{kn} & = & d_{kn} - d_{vn} & \mbox{for } 1 \leq k \leq n-1.\\ 
\end{array}
\]

It can be easily seen that $\mathbf{dy} \leq d_0$ is a convex
combination of  $\mathbf{ay} \leq a_0$ and $\mathbf{by} \leq b_0$:
\[(\mathbf{d},d_0)= \frac{d_{vn}}{d_{rs}+d_{vn}}(\mathbf{a},a_0) + \frac{d_{rs}}{d_{rs}+d_{vn}}(\mathbf{b},b_0).\]
Further, all three inequalities are pairwise nonequivalent; so it
remains to be shown that the inequalities $\mathbf{ay} \leq a_0$ and
$\mathbf{by} \leq b_0$ are valid for $\UPaththree$. 
This can be done by checking $a_{ij}+a_{jn} \leq a_0$ and
$b_{ij}+b_{jn} \leq b_0$ for all $1 \leq i,j \leq n-1$ with $i \neq j$.

Let $i$ and $j$ be distinct nodes in $\{1,\dots, n-1\}$. 

\vspace{0.5cm}
\noindent CASE 1: $[i,j], [j,n] \notin M$.

We have $a_{ij}=b_{ij}=d_{ij}$ and $a_{jn}=b_{jn}=d_{jn}$. Thus,
$a_{ij}+a_{jn} \leq a_0$ and $b_{ij}+b_{jn} \leq b_0$, since
$d_{ij}+d_{jn} \leq d_0$.

\vspace{0.5cm}
\noindent CASE 2: $[i,j] \in M$, $[j,n] \notin M$.

Since $0 < d_{ij} < d_0$, $d_{jn} \in \{0,d_0\}$, and $d_{ij}+d_{jn}
\leq d_0$, we deduce that $d_{jn}=0$. Hence, also $a_{jn}=b_{jn}=0$.
Since $a_{ij}=d_{ij}-d_{rs} < d_{ij}$, it follows that $a_{ij}+a_{jn}
\leq a_0$. Due to (a), $[i,n] \in M$, and since $d_{in} \geq d_{vn}$,
we deduce that $d_{ij} \leq d_0 - d_{vn}$. Thus,
$b_{ij}+b_{jn}=d_{ij}+d_{vn} \leq d_0=b_0$.

\vspace{0.5cm}
\noindent CASE 3: $[i,j] \notin M$, $[j,n] \in M$.

This implies that $a_{ij}=b_{ij}=d_{ij}=0$ and thus, $b_{ij}+b_{jn} \leq b_0$.
Due to (b), there is some internal node $l$ such that $[l,j] \in M$.
Since $d_{lj} \geq d_{rs}$, we deduce that $d_{jn} \leq d_0-d_{rs}$
and hence, $a_{ij}+a_{jn}=d_{jn}+d_{rs} \leq d_0=a_0$.

\vspace{0.5cm}
\noindent CASE 4: $[i,j],[j,n] \in M$.

Clear.

\vspace{0.5cm}

Thus, in all four cases, the inequalities  $\mathbf{ay} \leq a_0$ and
$\mathbf{by} \leq b_0$ are valid for $\UPaththree$. So we have shown that
$d_e \in \{0,d_0\}$ for all $e \in E$ and without loss of generality,
we may assume that $d_0=1$.

We resume: the facet defining inequality $\mathbf{dy} \leq d_0$ satisfies
(i)-(iii), $d_0=1$, and $d_e \in \{0,1\}$ for all $e \in E$. Note that
$d_{ln}=1$ for some internal node $l$ implies that $d_{il}=0$ for all
internal nodes $i \neq l$.

When $d_{in}=0$ for $1 \leq i \leq n-1$, we deduce that $d_e=1$ for
all internal edges $e \neq [u,w]$, i.e.,  $\mathbf{dy} \leq d_0$ is
equivalent to the nonnegativity constraint $y_{uw} \geq 0$.

When $d_{in}=1$ for all internal nodes $i \neq z$, we see that $d_e=0$
for all internal edges $e$. Then,  $\mathbf{dy} \leq d_0$ is
equivalent to the nonnegativity constraint $y_{zn} \geq 0$.

In all other cases, i.e., for $1 \leq \sum_{i=1}^{n-1} d_{in} \leq
n-2$, the inequality $\mathbf{dy} \leq d_0$ is not equivalent to a
nonnegativity constraint which implies that for each edge $e$ there is
a tight $[0,n]-3$-path containing $e$. Thus, $d_{ij}=1$ for all
internal edges $[i,j]$ for which $d_{in}=d_{jn}=0$. Therefore,
$\mathbf{dy} \leq d_0$ is a member of family (\ref{105}).
\end{proof}

Next, we turn to the polytopes $\UPath$ when $4 \leq p \leq n-1$.
The integer points in $\UPath$ are characterized by the following
model:
\begin{align}
y_{0n} & =0 \label{uzero}\\
y(\delta(0)) & =1 \label{d0} \\
y(\delta(n)) & =1 \label{dn} \\
y(\delta(j)) & \leq 2 & \forall \,\,\,j \in V\setminus \{0,n\}
\label{udegree}\\
y(\delta(j) \setminus \{e\}) - y_e & \geq 0 & \forall \,\,\, j \in V
\setminus \{0,n\}, e \in \delta(j), \label{parity} \\
y((S:V \setminus S)) & \geq y(\delta(j)) & \forall \,\,\, S \subset V,
3 \leq |S| \leq n-2, \label{uosmincut}\\
\nonumber & &  0,n \in S, j \in V \setminus S\\
y(E) &=p \label{ucard}\\
x_{e}  \in \{0,1\} & & \forall \,\,\, e \in E. \label{uinteger}
\end{align}
Here, for any node sets $S,T$ of
$V$, $y((S:T))$ is short for $\sum_{i \in S} \sum_{j \in T} y_{ij}$,
where the summation does not extend over loops $(i,i)$ for $i \in S
\cap T$.

The \emph{parity} constraints (\ref{parity}) together with the
degree (\ref{udegree}) and the
integrality constraints (\ref{uinteger}) ensure that every internal
node has degree 0 or 2. Hence, constraints (\ref{uzero}) - (\ref{parity})
and the integrality constraint (\ref{uinteger}) are satisfied by the
incidence vector of the node disjoint union of a simple $[0,n]$-path and
simple cycles on the set of internal nodes. The one-sided min-cut
inequality
(\ref{uosmincut}) is satisfied by the incidence vectors of simple
$[0,n]$-paths but violated by the incidence vectors of the union of a simple $[0,n]$-path and
simple cycles. Finally, the cardinality constraint (\ref{ucard})
excludes all incidence vectors of $[0,n]$-paths which have a length
that is not equal to $p$. 

\begin{Lemma} \label{U1}
Let $4 \leq p \leq n-1$ and $n \geq 6$. If the equation
\[\mathbf{cy}=c_0\]
is satisfied by all $[0,n]-p$-paths, then there are $\alpha,
\beta, \gamma$, such that $c_{0i}=\alpha$, $c_{in}=\beta$ for all $i
\in \{1,\dots,n-1\}$ and $c_{ij}=\gamma$ for all $i,j \in \{1,\dots,n-1\}$.
\end{Lemma}

\begin{proof}
Set $S:=\{1,\dots,n-1\}$ and let $i,j,k,l$ be any distinct nodes in $S$ and consider
any $[0,n]-p$-path $P$ that uses the edges $[i,j],[j,k]$ but does not
visit node $l$. Replacing node $j$ by node $l$ yields $c_{ij}+c_{jk}=
c_{il}+c_{kl}$. Next, consider any  $[0,n]-p$-path $P'$ that uses the
edges $[j,i],[i,l]$ but does not visit the node $k$. Replacing node
$i$ by node $k$ yields $c_{ij}+c_{il}=c_{jk}+c_{kl}$. We deduce that
$c_{ij}=_{kl}$ and since $|S| \geq 5$, we see that
$c_{ij}=c_{kl}$ for all distinct nodes $i,j,k,l \in S$.
Denoting this common value by $\gamma$, it follows
immediately that there are $\alpha, \beta$ with $c_{0i}=\alpha$ and
$c_{in}=\beta$ for all $i \in S$.
\end{proof}

We are now well prepared to determine the dimension of $\UPath$
depending on $n$ and $p$.
For the sake of completeness we determine also the dimension of
$\UPath$ when $p=n$.
 
\begin{Theorem}
Let $n \geq p \geq 4$. Then
\[
\dim \UPath = 
\left \{ 
\begin{array}{rl}
|E| -4 & \mbox{if } p \leq n-1,\\
|E|-n-2  & \mbox{if } p=n \geq 4.\\
\end{array}
\right .
\]
\end{Theorem}

\begin{proof}
Using a convex hull code we see
that $\dim \UPathff = 11$. Next, suppose that $n \geq 6$ and $4
\leq p \leq n-1$. We will show that (\ref{uzero})-(\ref{dn}) and
(\ref{ucard}) is a minimal equality subsystem for $\UPath$. Since
the equations (\ref{uzero})-(\ref{dn}) and
(\ref{ucard}) are linearly independent, $\dim \UPath \leq
\frac{(n+1)n}{2}-4$. It remains to be shown that any equation that is
satisfied by all $\mathbf{y} \in \UPath$ is a linear combination of
(\ref{uzero})-(\ref{dn}) and (\ref{ucard}). Let $\mathbf{cy}=c_0$ be
such an equation. By Lemma \ref{U1}, there are $\alpha, \beta, \gamma$
 with $c_{0i}=\alpha$, $c_{in}=\beta$ for all internal nodes $i$ and
 $c_{ij}=\gamma$ for all internal nodes $i \neq j$. Thus,
\[\begin{array}{rcl}
(\mathbf{cy},c_0) & = & \gamma ( y(E), p) \\
 & & + (\alpha - \gamma) (y(\delta(0)), 1)\\
& & + (\beta-\gamma) ( y(\delta(n)), 1) \\
& & + (c_{0n}+\gamma-\alpha-\beta)  (y_{0n},0).
\end{array}
\]

Finally, let $p=n \geq 4$. Theorem 7 of Gr\"otschel and Padberg
\cite{GP} implies that the dimension of the traveling salesman
polytope $Q_T^{n+1}$ defined on the complete graph on node set $V$ is equal
to $|E|-n-1$ for $n \geq 2$ and Theorem 8 of the same authors
\cite{GP} says that the inequalities $x_{e} \leq 1$ induce facets
$F_e$ of $Q_T^{n+1}$ for $n \geq 3$. Since $F_{0n}$ is isomorphic to
$\UPathn$, we obtain the required result. 
\end{proof}

A valid inequality
$\mathbf{cx} \leq c_0$
for the \path polytope $\Path$ is said to be \symmetric if
$c_{ij}=c_{ji}$ for all $1 \leq i < j \leq n-1$. It is easy to see
that the undirected counterpart $\mathbf{\bar{c}y} \leq c_0$ of a
\symmetric inequality $\mathbf{cx} \leq c_0$
(obtained by setting $\bar{c}_{0i}=c_{0i},\bar{c}_{in}=c_{in}$ for all
internal nodes $i$ and $\bar{c}_{ij}=c_{ij}=c_{ij}$ for all $1 \leq i
< j \leq n-1$) is facet defining for $\UPath$ if $\mathbf{cx} \leq
c_0$ is facet defining for $\Path$ (cf. \cite{HO}). The argument that
can be used to prove the statement is the following: assuming that 
 $\mathbf{\bar{c}y} \leq c_0$ does not induce a facet of $\UPath$,
 then there is a facet inducing inequality  $\mathbf{\bar{d}y} \leq d_0$ for $\UPath$ such that
 $\{\mathbf{y} \in \UPath | \mathbf{\bar{c}y} = c_0 \} \subsetneq
\{\mathbf{y} \in \UPath | \mathbf{\bar{d}y} = d_0 \} $. But then
 $\{\mathbf{x} \in \Path | \mathbf{cy} = c_0 \} \subsetneq
\{\mathbf{x} \in \Path | \mathbf{dy} = d_0 \} $, where
$\mathbf{dx} \leq d_0$ is the directed counterpart of
$\mathbf{\bar{d}y} \leq d_0$ (obtained by setting $d_{0i}=
\bar{d}_{0i}$, $d_{in} = \bar{d}_{in}$ for all $i \in \{1,\dots,n-1\}$
and $d_{ij}=d{ji}=\bar{d}_{ij}$ for all $1 \leq i < j \leq n-1$).

Since the degree constraint
(\ref{degree1}) and the cut inequalities (\ref{mincut}), (\ref{osmincut1}),
(\ref{oddmaxcut}), and (\ref{evenmaxcut}) are \symmetric, their
undirected counterparts are facet defining for $\UPath$.

\begin{Corollary}
Let $4 \leq p <n$.
\begin{itemize}
\item[(i)] The degree constraint $y(\delta(j)) \leq 2$ induces a facet
  of $\UPath$ for every internal node $j$ of $G$.
\item[(ii)] Let $S \subset V$ and $0,n \in S$. The min-cut inequality $y((S:V \setminus S)) \geq 2$
  induces a facet of $\UPath$ if $3 \leq |S|
  \leq p$.
\item[(iii)]  Let $S \subset V$ and $0,n \in S$. The one-sided min-cut inequality $y((S:V \setminus S))
  \geq y(\delta(j))$ defines a facet of $\UPath$ for every node $j \in
  V \setminus S$.
\item[(iv)]   Let $S \subset V$ and $0,n \in S$. The max-cut
  inequality $y((S:T)) \leq p-1 $ defines a facet of $\UPath$ if
$p$ is odd, $S \setminus \{n\}> p/2$, and $T > p/2$.
\item[(v)]  Let $S \subset V$ and $0 \in S$ and $n \in T$. The
  max-cut inequality $y((S:T)) \leq p/2$ induces a facet of $\UPath$
  if $p$ is even, $|S| > p/2$, and $|T|> p/2$. 
\end{itemize}
\hfill $\Box$
\end{Corollary}

Finally, we show that the nonnegativity constraints $x_e \geq 0$
define facets of the $[0,n]-p$-path polytope $\UPath$.

\begin{Theorem}
Let $4 \leq p < n$. The nonnegativity constraint
\begin{equation} \label{unon}
y_e \geq 0
\end{equation}
defines a facet of the $[0,n]-p$-path polytope $\UPath$ for all edges
$e \neq [0,n]$ of $K_{n+1}$.
\end{Theorem}

\begin{proof}
When $n \leq 5$, (\ref{unon}) can be seen to be facet defining using
a convex hull code; so assume that $n \geq 6$. Let $\mathbf{cy}=c_0$
be an equation 
that is satisfied by every $\mathbf{y} \in \UPath$ with $y_{e}=0$.
Since the lineality space of $\UPath$ is determined by the equations
(\ref{uzero})-(\ref{dn}) and (\ref{ucard}), we may assume that
$c_{0n}=0$, $c_{0m}=c_{mn}=0$ for some internal node $m$ with $[0,m]
\neq e \neq [m,n]$, and $c_f=0$ for some internal edge $f \neq e$.

Let $g=[i,j], h=[k,l] \in E
\setminus \{e\}$ be not adjacent edges. Without loss of generality,
we may assume that the nodes $j$ and $l$ are not incident with edge
$e$. 
Let $P$ be any tight
$[0,n]-p$-path that uses the edges $[i,j], [j,k]$ but does not visit
node $l$. Replacing node $j$ by node $l$ yields another tight path and
hence, $c_{ij}+c_{jk}=c_{il}+c_{lk}$. Next, consider any tight
$[0,n]-p$-path $P'$ that uses the edges $[j,i],[i,l]$ and does not visit
node $k$. Replacing node $i$ by node $k$ yields another tight path and
thus, $c_{ij}+c_{jk}=c_{il}+c_{lk}$. Adding both equations, we obtain
$c_g=c_h$, and since $|V \setminus \{0,n\}| \geq
5$, this implies $c_g=c_h$ for all internal edges $g,h$ that are not
equal to $e$. Now it is easy to see that also $c_{0i}=c_{0j}$ and
$c_{kn}=c_{ln}$ for all edges $[0i],[0j],[kn],[ln]$ not equal to $e$.
Since  $c_{0m}=c_{mn}=0$ and $c_f=0$,
it follows that $c_g=0$ for all edges $g \neq e$ which implies also
$c_0=0$.
Hence, $\mathbf{cx}=c_0$ is simply $c_e y_e=0$.
\end{proof}

\subsection{Facets of the lower and upper directed \path polytopes} \label{relax}

\begin{Theorem}[cf. Theorem 18 of Hartmann and \oez \cite{HO}]
Let $\mathbf{cx} \leq c_0$ induce a facet of the \path polytope $\Path$,
where $4 \leq p < n$. If $\mu$ is the smallest (largest) value such that
\begin{equation}
\mu x(A) + \mathbf{cx} \leq \mu p + c_0 \label{relaxlift}
\end{equation}
is valid for the lower (upper) \path polytope, then \emph{(\ref{relaxlift})}
is facet inducing for the lower (upper) \path polytope.
\hfill $\Box$
\end{Theorem}

\begin{Corollary}
Let $4 \leq p < n$.
The nonnegativity constraints \emph{(\ref{nn})}, 
degree constraints  \emph{(\ref{degree1})},
one-sided min-cut inequalities  \emph{(\ref{osmincut1})}, 
 max-cut inequalities  \emph{(\ref{Roddmaxcut}) -  (\ref{Revenmaxcut2})},
jump inequalities  \emph{(\ref{jump})},
and cardinality-path inequalities  \emph{(\ref{cardpath})} 
are facet defining for the lower \path polytope, if the accordant
conditions hold. 

\hspace{3cm}
\hfill $\Box$
\end{Corollary}

\begin{Corollary}
Let $4 \leq p < n$, $S \subset V$, and $0,n \in S$.
The inequality
\begin{equation}
x(A) - x((S:V \setminus S)) \leq p-1 \label{mincutl}
\end{equation}
induces a facet of the lower \path polytope if and only if 
$|S| \leq p$ and $|V \setminus S| \geq 2$.
\end{Corollary}

\begin{proof}
The inequality  (\ref{mincutl}) is derived from the min-cut inequality
  (\ref{mincut}) with parameter $\mu=-1$. Hence it is facet defining, if 
$3 \leq |S| \leq p$ and $|V \setminus S| \geq 2$.
When $S=\{0,n\}$, (\ref{mincutl}) is equivalent to the 
cardinality constraint $x(A) \leq p$ and hence facet defining for the lower \path polytope. 

Conversely, when $|S| \geq p+1$,  (\ref{mincutl}) is no longer valid,
and when $|V \setminus S|=1$, $n \leq p$, a contradiction.
\end{proof}

\newpage
\noindent
\small{E-mail address: \texttt{stephan@math.tu-berlin.de}}

\vspace{0.5cm}
\noindent
TU Berlin, MA 3-1, Stra{\ss}e des 17. Juni 136, 10623 Berlin, Germany


\begin{thebibliography}{99}
\bibitem[1]{BO} E. Balas and M. Oosten, \emph{On the cycle polytope of a directed graph}, Networks 36, No. 1 (2000), pp. 34-46.
\bibitem[2]{Bauer} P. Bauer, \emph{A Polyhedral Approach to the Weighted Girth Problem}, Aachen 1995.
\bibitem[3]{BLS} P. Bauer, J.T. Linderoth, and M.W.P. Savelsbergh, \emph{A branch and cut approach to the cardinality constrained circuit problem}, Mathematical Programming, Ser. A 91 (2002), pp. 307-348.
\bibitem[4]{Christof} T. Christof, \emph{Ein Verfahren zur Transformation zwischen Polyederdarstellungen}, Diplomarbeit, Universit\"at Augsburg, 1991.
\bibitem[5]{CP} C. Coullard and W.R. Pulleyblank, \emph{On cycle cones and polyhedra}, Linear Algebra Appl. 114/115 (1989), pp. 613-640.
\bibitem[6]{Dahl} G. Dahl and L. Gouveia, \emph{On the directed hop-constrained shortest path problem}, Operations Research Letters 32 (2004),
pp. 15-22.
\bibitem[7]{DR} G. Dahl and B. Realfsen, \emph{The
    Cardinality-Constrained Shortest Path Problem in 2-Graphs},
  Networks 36 No. 1 (2000), pp. 1-8.
\bibitem[8]{polymake} E. Gawrilow and M. Joswig, \emph{polymake: A
    framework for analyzing convex polytopes}. In: G. Kalai and G.M. Ziegler (eds.): Polytopes ~ Combinatorics and Computation
(DMV-Seminars, pp. 43~74) Basel: Birkh¨auser-Verlag Basel 2000,
see also http://www.math.tu-berlin.de/polymake
\bibitem[9]{Groetschel} M. Gr\"otschel, \emph{Polyedrische Charakterisierungen kombinatorischer Optimierungsprobleme},
Mathematical Systems in Economics 36, Meisenheim am Glan, 1977.
\bibitem[10]{HO} M. Hartmann and \"O. \oez, \emph{Facets of the $p$-cycle polytope}, Discrete Applied Mathematics 112 (2001), pp. 147-178.
\bibitem[11]{GP} M. Gr\"otschel and M.W. Padberg, \emph{Polyhedral theory}, in: E.L. Lawler et al (eds.),
\emph{The traveling salesman problem. A guided tour of combinatorial optimization}, Chichester, New York, and others, 1985, pp. 
251-305.
\bibitem[12]{MN1} J. Maurras and V.H. Nguyen, \emph{On the linear description of the 3-cycle polytope}, European Journal of Operational Research, 1998.
\bibitem[13]{MN2} J. Maurras and V.H. Nguyen, \emph{On the linear
    description of the k-cycle polytope, $PC_n^k$}, International
  Transactions in Operational Research 8 (2001), pp. 673-692.
\bibitem[14]{NW}G.L. Nemhauser and L.A. Wolsey, Integer and Combinatorial Optimization, Wiley: New York, 1988.
\bibitem[15]{Schrijver2003} A. Schrijver, \emph{Combinatorial
    Optimization}, Vol. A, Berlin et al., 2003.
\bibitem[16]{Stephan} R. Stephan, \emph{Polytopes associated with length restricted directed circuits}, Master Thesis, Technische Universit\"at Berlin, 2005.
\end{thebibliography}
\end{document}